\documentclass[final]{siamart220329}
\usepackage{amssymb, amsfonts, bm, macrodefs, algorithmic, graphicx, subcaption, nicematrix, enumitem}

\title{Collect, Commit, Expand: Efficient CPQR-Based Column Selection for Extremely Wide Matrices}

\author{Robin Armstrong\thanks{Cornell University, Center for Applied Mathematics, Ithaca, NY (\email{rja243@cornell.edu}).} \and Anil Damle\thanks{Cornell University, Department of Computer Science, Ithaca, NY (\email{damle@cornell.edu}).}}

\newcommand{\delete}[1]{}

\begin{document}

\maketitle

\begin{abstract}
Column-pivoted QR (CPQR) factorization is a computational primitive used in numerous applications that require selecting a small set of ``representative'' columns from a much larger matrix. These include applications in spectral clustering, model-order reduction, low-rank approximation, and computational quantum chemistry, where the matrix being factorized has a moderate number of rows but an extremely large number of columns. We describe a modification of the Golub-Businger algorithm which, for many matrices of this type, can perform CPQR-based column selection much more efficiently. This algorithm, which we call CCEQR, is based on a three-step ``collect, commit, expand'' strategy that limits the number of columns being manipulated, while also transferring more computational effort from level-2 BLAS to level-3. Unlike most CPQR algorithms that exploit level-3 BLAS, CCEQR is deterministic, and provably recovers a column permutation equivalent to the one computed by the Golub-Businger algorithm. Tests on spectral clustering and Wannier basis localization problems demonstrate that on appropriately structured problems, CCEQR can significantly outperform GEQP3.
\end{abstract}


\section{Introduction} The column subset selection problem (CSSP) appears in a remarkably wide range of applications. For example, certain nonlinear model order reduction techniques involve finding a small set of ``informative'' state components that capture the nonlinear term \cite{chaturantabut_deim}. Some spectral clustering algorithms involve finding a small subset of data points containing one representative from each cluster \cite{damle_robust_clustering}. And many low-rank approximation techniques require finding a small column subset whose span approximates the range of a larger matrix \cite{murray_rnla_survey}. All of these problems can be treated as instances of CSSP.

Given a matrix $\mA \in \C^{m \times n}$ and an integer $1 \leq k \leq \min\{ m,\, n \}$, CSSP entails finding $k$ columns, called ``skeleton columns,'' that are maximally large and linearly independent. Skeleton column selection can be formalized as finding $\vs \in [n]^k$ that maximizes $\sigma_\mathrm{min}(\mA(:,\vs))$, which is an NP-hard combinatorial optimization problem \cite{civril_2009}. In the communities of numerical linear algebra and theoretical computer science, significant effort has been devoted to finding efficient CSSP approximation algorithms over the course of several decades. Most of the resulting algorithms select columns using random sampling \cite{deshpande_rademacher_volumesampling, deshpande_rademacher_projectiveclustering, mahoney_drineas_muthukrishnan, frieze_monte_carlo, mahoney_cur} or pivoted matrix factorizations \cite{businger_linear_least_squares, chandrasekaran_1994, drmac_qdeim, golub_klema_stewart, gu_srrqr, hong_pan, sorensen_embree_deim}, which themselves may be randomized \cite{rgks, martinsson_rokhlin_tygert_randid, martinsson_randomized_qrcp, saibaba_rdeim, woolfe_srft}. The Golub-Businger algorithm for column-pivoted QR factorization \cite{businger_linear_least_squares} is a widely used column selection tool, owing to the fact that high-performance implementations exist and can be easily called from the standard libraries of most scientific computing languages. Although it has poor worst-case performance \cite{gu_srrqr, kahan_nla}, the empirical performance of this algorithm on realistic problems is usually quite strong \cite{damle_rank_revealing}.

A major difficulty in the Golub-Businger algorithm is that in a na\"ive implementation, Householder reflectors must be interwoven with column pivots. This precludes blocking the reflectors into compact forms that would allow a group of reflectors to be applied simultaneously with level-3 BLAS routines \cite{bischof_van_loan_wy, schreiber_van_loan_compactwy}, as opposed to sequentially with much slower BLAS-2 operations. Like many authors before us, we will demonstrate a modification of the Golub-Businger algorithm that reduces BLAS-2 computation by allowing blocked application of Householder reflectors. Most prior work has achieved this by limiting BLAS-2 work to a small number of rows, either by transferring the entire pivot calculation onto a sketched matrix with shorter columns \cite{martinsson_rokhlin_tygert_randid, woolfe_srft}, by reflecting single rows at a time to determine pivots while accumulating reflectors into compact WY blocks \cite{orti_blas3_cpqr}, or by generating reflectors and pivots in blocks using repeated sketches of the trailing submatrix \cite{martinsson_randomized_qrcp}. Given that the Golub-Businger algorithm has $\cO(m^2 n)$ complexity for an $m \times n$ matrix, these techniques can dramatically speed up computation when $m$ is moderate or large compared to $n$.

This paper, however, is concerned with situations where $n$ vastly exceeds $m$. Several applications give rise to this sort of column subset selection problem. For example, spectral clustering \cite{damle_robust_clustering} generates CSSP instances where each row corresponds to a cluster and each column corresponds to a single data point. Model order reduction \cite{chaturantabut_deim} generates CSSP instances where each row corresponds to a coordinate of the reduced model and each column corresponds to a coordinate of the full model. In these examples and more, the number of columns may exceed the number of rows by several orders of magnitude. Unless the desired number of skeleton columns is much smaller than the number of rows, strategies based limiting BLAS-2 effort to a small number of rows do not effectively address the main bottleneck for these problems.

We will address this bottleneck using an algorithm that we call CCEQR (``Collect, Commit, Expand'' QR). This algorithm repeatedly ``collects'' a set of heuristically good pivot candidates, ``commits'' a provably good subset of the candidates into the pivot sequence, and then ``expands'' the candidate search to a wider range of columns. Unlike column selection strategies for wide matrices that use distributed parallelism \cite{bishof_local_pivoting, demmel_tournament_pivoting}, CCEQR runs on a single machine, and unlike methods based on random sampling \cite{deshpande_rademacher_volumesampling, deshpande_rademacher_projectiveclustering, mahoney_drineas_muthukrishnan, frieze_monte_carlo, mahoney_cur}, it is fully deterministic. CCEQR is distinguished from both classes of methods mentioned above by the fact that it provably selects the same column subset as the Golub-Businger algorithm (cf.\ \cref{theorem:cceqr_equivalence}). In terms of efficiency, the major advantages of CCEQR are that (1) Householder reflectors are only applied with BLAS-2 during the ``commit'' step, where they are limited to the small set of candidate columns, and (2) Householder reflections in the ``expand'' step (performed with BLAS-3) are only applied to a small subset of ``tracked'' columns. We explain the details of this algorithm in \cref{section:cceqr}. In \cref{section:experiments} we will show that for an appropriately structured problem, CCEQR can perform as much as an order of magnitude faster than the standard implementation of the Golub-Businger algorithm in LAPACK.

By an ``appropriately structured problem,'' we mean two things. First, we are referring to problems where the relevant information to compute is the column permutation, rather than the full column-pivoted QR factorization. This includes the aforementioned problems in model order reduction \cite{chaturantabut_deim}, spectral clustering \cite{damle_robust_clustering}, and low-rank approximation \cite{murray_rnla_survey}, in addition to many others. Although CCEQR can be easily modified to compute a full CPQR factorization $\mA\mPi = \mQ\mR$, much of its efficiency derives from the fact that the main algorithm computes only $\mQ$, $\mPi$, and \emph{a few} columns of $\mR$. Second, matrices that are ``appropriately structured'' for CCEQR have the property that their column norm distribution is mostly concentrated on a small subset of column indices. This assumption is both reasonable and interpretable in many application areas. For instance, in spectral clustering \cite{damle_robust_clustering}, column norms in the relevant matrix are determined by the likelihoods of the corresponding data points under various clusters \cite{schiebinger_clustering_geometry}. For these problems, CCEQR will perform well when a large number of data points are far from the center of any cluster. Density functional theory generates a column subset selection problem where column norms in the relevant matrix are directly related to electron densities over a discretized spatial grid \cite{damle_scdm, lin_lu_est}. These norms will decay rapidly under mild physical assumptions about the chemical system \cite{benzi_projector_decay}. In \cref{section:experiments}, we will empirically demonstrate how variations in problem structure impact the performance of CCEQR. In the worst cases, the runtime of CCEQR will still be comparable to (albeit slower than) the LAPACK implementation of the Golub-Businger algorithm.

\subsection{Code Availability}
A reference implementation of CCEQR is available at \url{https://github.com/robin-armstrong/CCEQR.jl/releases/tag/v1.1.0}. Code to reproduce the numerical experiments in this paper can found at \url{https://github.com/robin-armstrong/cceqr-experiments/releases/tag/arXiv-v2}.

\subsection{Notation} Matrices will be denoted with bold capital letters, vectors with bold lowercase letters, and all other objects with un-bolded letters. We use \textsc{Matlab} notation to slice vectors and matrices by row or column. Thus, if $\vx$ is a vector, then $\vx(j)$ is its $j\nth$ entry. If $\mA$ is a matrix, then $\mA(i, :)$ denotes its $i\nth$ row, $\mA(:, j)$ its $j\nth$ column, and $\mA(i,\, j)$ its $(i,\, j)$ element. When $\vj = [j_1,\, \ldots,\, j_r]$ is a vector of column indices, we will write $\vx(\vj)$ to denote the vector $[\vx(j_1),\, \ldots,\, \vx(j_r)]$, $\mA(:,\vj)$ to denote the matrix $[\mA(:, j_1) \: \cdots \: \mA(:, j_r)]$, and similar notation applies to rows. The symbol $[n]$ represents the set $\{ 1,\, 2,\, \ldots,\, n \}$, and when $p \leq q$ are integers, the symbol $p \mcol q$ represents the vector $[p,\, p + 1,\, \ldots,\, q - 1,\, q]$.

Indexed matrix products will be ordered from left to right. Thus, if $\mA_1,\, \ldots,\, \mA_k$ are matrices of conformal dimensions, then $\prod_{i = 1}^k \mA_i \defeq \mA_1\mA_2 \ldots \mA_k$.
\section{Background} \label{section:background}
Here we will review the necessary background material for this paper, including column-pivoted and rank-revealing QR factorizations, the Golub-Businger algorithm, blocking techniques for Householder reflectors, and accelerated column selection strategies based on randomization and parallelization.

\subsection{CPQR Factorizations} Given $\mA \in \C^{m \times n}$ and $k \in \{ 1,\, \ldots,\, \min\{ m,\, n \} \}$, a column-pivoted QR (CPQR) factorization of $\mA$ has the form
\begin{equation*}
\mQ^*\mA\mPi = \begin{bmatrix}
\mR_{11} & \mR_{12} \\
\mZero & \mR_{22}
\end{bmatrix},
\end{equation*}
where $\mQ \in \C^{m \times m}$ is unitary, $\mPi \in \{ 0,\, 1 \}^{n \times n}$ is a permutation, and $\mR_{11} \in \C^{k \times k}$ is upper-triangular. The first $k$ columns of $\mPi$ define an index vector $\vs$ such that $\mA\mPi(\mcol,1 \mcol k)= \mA(\mcol,\vs)$. CPQR algorithms construct $\mPi$ in order to satisfy two basic conditions, which can be informally stated as (1) the columns of $\mA(\mcol,\vs)$, called the ``skeleton columns,'' are large and well-conditioned, and (2) $\range \mA(\mcol,\vs)$ is a low-dimensional approximation of $\range \mA$. Note that $\range \mA(\mcol,\vs) = \range \mQ(\mcol, 1 \mcol k)$, and that the columns of $\mR_{11}$ are the coordinates of the skeleton columns in the basis $\mQ(\mcol, 1 \mcol k)$. The non-skeleton columns can be orthogonally decomposed into their components in $\range \mQ(\mcol, 1 \mcol k)$, whose coordinates are the columns of $\mR_{12}$, and their components in $(\range \mQ(\mcol, 1 \mcol k))^\perp$, whose coordinates are the columns of $\mR_{22}$.

More formally, CPQR algorithms choose $\mPi$ such that
\begin{align}
&\sigma_\mathrm{min}(\mA(:,\vs)) = \sigma_\mathrm{min}(\mR_{11}) \geq \frac{\sigma_k(\mA)}{q(n,\, k)} \quad\text{and} \label{eqn:rrqr_factorization_1} \\
&\| \mA - \mA(:,\vs)\mA(:,\vs)\pinv\mA \|_2 = \sigma_\mathrm{max}(\mR_{22}) \leq q(n,\, k)\sigma_{k + 1}(\mA), \label{eqn:rrqr_factorization_2}
\end{align}
where $q$ is a function determined by the algorithm's design that bounds the suboptimality of the factorization\footnote{Here and throughout this paper, $\sigma_i(\cdot)$ denotes the $i\nth$ largest singular value of a matrix, $\sigma_\mathrm{max}(\cdot) = \sigma_1(\cdot)$, and $\sigma_\mathrm{min}(\cdot)$ denotes the $\min\{ m,n \}\nth$ largest singular value of a matrix with $m$ rows and $n$ columns.}$^,$\footnote{This function bounds suboptimality because $\sigma_\mathrm{min}(\mA(:,\vj)) \leq \sigma_k(\mA)$ for any $\vj \in [n]^k$ (a consequence of singular value interlacing), and $\| \mA - \mA(:,\vj)\mB \|_2 \geq \sigma_{k + 1}(\mA)$ for any $(\vj,\, \mB) \in [n]^k \times \C^{k \times n}$ (a consequence of the Eckart-Young theorem).}. If $q$ is bounded by a low-degree bivariate polynomial then the algorithm is said to compute a ``rank-revealing'' QR factorization \cite{chandrasekaran_1994, damle_rank_revealing, gu_srrqr, hong_pan}. Although they have a strong near-optimality guarantee, rank-revealing QR factorization algorithms are complex to implement and, crucially, are not implemented in LAPACK. This paper focuses on the much more commonly used Golub-Businger algorithm \cite{businger_linear_least_squares}, implemented in LAPACK as GEQP3.

\subsection{The Golub-Businger Algorithm} Rather than maximizing $\sigma_\mathrm{min}(\mA(:,\vs))$ over $\vs \in [n]^k$, the Golub-Businger algorithm \cite{businger_linear_least_squares} constructs a column permutation which solves a related but more computationally tractable problem: that of reducing a matrix to $k$-Golub-Businger form.
\begin{definition} \label{def:gb_form}
Let $\mR \in \C^{m \times n}$ and $k \leq \min\{ m,\, n \}$. We say that $\mR$ has \emph{$k$-Golub-Businger form}, or \emph{$\mathrm{GB}(k)$ form}, if its first $k$ columns are upper-triangular and satisfy
\begin{equation}
|\mR(i,i)| = \max_{i \leq j \leq n} \| \mR(i \mcol m, j) \|_2 \label{eqn:gb_requirement}
\end{equation}
for $1 \leq i \leq k$. We observe the convention that all matrices are in $\mathrm{GB}(0)$ form.
\end{definition}

The first $k$ columns of a $\mathrm{GB}(k)$-form matrix represent a ``greedily optimal'' column subset for large norms and good conditioning, in the following sense: for $1 \leq i \leq k$, the $i\nth$ column has maximal norm in the subspace orthogonal to the previous $i - 1$. Given $\mA \in \C^{m \times n}$, the Golub-Businger algorithm computes a unitary $\mQ$ and a permutation $\mPi$ such that $\mR \defeq \mQ^*\mA\mPi$ has GB$(k)$ form. This is accomplished by setting $\mR \leftarrow \mA$, $\mQ \leftarrow \mI_{m \times m}$, $\mPi \leftarrow \mI_{n \times n}$, and then repeating the following actions for $i = 1,\, \ldots,\, k$:
\begin{enumerate}
\item Locate the column with largest norm orthogonal to $\range \mR(i \mcol m, 1 \mcol (i-1))$. Due to the upper-triangular structure in $\mR(i \mcol m, 1 \mcol (i-1))$, this corresponds to finding $j_\mathrm{max} \in \argmax_{i \leq j \leq n} \| \mR(i \mcol m,j) \|_2$.
\item Move this column to index $i$, and record this permutation by modifying $\mPi$.
\item Apply a rotation that creates zeros in the $((i+1) \mcol m,i)$ block of $\mR$, and record this rotation by modifying $\mQ$.
\end{enumerate}
The LAPACK routine GEQP3 implements this algorithm with $k = \min\{m,\, n\}$.

We summarize the Golub-Businger algorithm in \cref{alg:golub_businger}. The efficiency and numerical stability of this algorithm depends critically on the use of Householder reflections in lines \eqref{line:GEQP3_householder} through \eqref{line:GEQP3_Q_update}, the details of which are discussed in \cref{subsec:householder_review}. Fast implementations of this algorithm such as GEQP3 include several performance optimizations which the pseodocode in \cref{alg:golub_businger} does not capture. For example, $\mQ$ is represented in compressed form in terms of its Householder vectors, $\mPi$ is stored as a vector, and recursive update formulas are used to avoid recalculating column norms from scratch at line \eqref{line:GEQP3_argmax}. See \cite[sec 5.4.2]{golub_van_loan} for details, and \cite{zlatko_bujanovic_rrqr_failure} for a discussion of numerical stability in the column-norm updates.
\begin{algorithm}
\caption{\textsc{Golub-Businger Algorithm}} \label{alg:golub_businger}
\begin{algorithmic}[1]
\STATE \textbf{inputs:} $\mA \in \C^{m \times n}$ (input matrix), $k \leq \min\{ m,\, n \}$ (desired skeleton size).
\STATE \textbf{outputs:} $\mQ,\, \mR,\, \mPi$ (factors of a CPQR factorization with $\mR$ in GB$(k)$ form).
\item[]
\STATE initialize $\mQ \leftarrow \mI_{m \times m},\, \mR \leftarrow \mA$, and $\mPi \leftarrow \mI_{n \times n}$.
\FOR{$i = 1,\, 2,\, \ldots,\, k$}
\STATE find $j_\mathrm{max} \in \argmax_{j \geq i} \| \mR(i \mcol m,j) \|_2$. \label{line:GEQP3_argmax}
\STATE swap columns $i$ and $j_\mathrm{max}$ in $\mR$.
\STATE swap columns $i$ and $j_\mathrm{max}$ in $\mPi$.
\IF{$i < k$}
\STATE compute $\vv_i,\, \tau_i \leftarrow \mathsf{Householder}(\mR(:,i), i)$. \label{line:GEQP3_householder}
\STATE update $\mR \leftarrow \mR - \tau_i\vv_i(\vv_i^*\mR)$. \label{line:GEQP3_R_update}
\STATE update $\mQ \leftarrow \mQ - \tau_i (\mQ\vv_i)\vv_i^*$. \label{line:GEQP3_Q_update}
\ENDIF
\ENDFOR
\RETURN $\mQ,\, \mR,\, \mPi$.
\end{algorithmic}
\end{algorithm}

The Golub-Businger algorithm satisfies \cref{eqn:rrqr_factorization_1,eqn:rrqr_factorization_2} with $q(n, k) = 2^k \sqrt{n - k}$ \cite{gu_srrqr}, though the actual performance of this algorithm is usually much better than the exponential factor in $k$ would suggest \cite{damle_rank_revealing}. This exponential factor reflects the fact that for certain pathological matrices, notably the $n \times n$ Kahan matrix with $k = n - 1$, the Golub-Businger algorithm selects a column subset with highly suboptimal conditioning \cite{kahan_nla}. Such matrices are almost never encountered in real-world applications.

\subsection{Householder Reflections and Compact WY Form} \label{subsec:householder_review} The unitary matrix $\mQ \in \C^{m \times m}$ produced by \cref{alg:golub_businger} has the form $\mQ = \mQ_1\mQ_2 \ldots \mQ_k$, where $\mQ_i \in \C^{m \times m}$ is a unitary matrix such that $\mQ_i^*$ creates zeros in the $((i+1) \mcol m,i)$ block of $\mR$. Efficient implementations of QR factorization algorithms do not form $\mQ_i$ explicitly, but instead represent it in terms of a $\vv_i \in \C^m$ and $\tau_i \geq 0$ such that $\mQ_i = \mI - \tau_i\vv_i\vv_i^*$. In \cref{alg:golub_businger}, $\vv_i$ and $\tau_i$ are computed at line \eqref{line:GEQP3_householder}, where $\vv,\, \tau \leftarrow \mathsf{Householder}(\vx,\, i)$ produces a ``Householder vector'' $\vv \in \C^m$ and a ``Householder scalar'' $\tau \geq 0$ such that
\begin{align*}
&\vv(1) = \ldots = \vv(i - 1) = 0, \\
&\vv(i) = 1, \\
&(\mI - \tau\vv\vv^*)^2 = \mI, \\
&(\mI - \tau\vv\vv^*)\vx = [\vx(1),\, \ldots,\, \vx(i-1),\, \mu,\, 0,\, \ldots,\, 0],\quad |\mu| = \pm\| \vx(i \mcol m) \|_2.
\end{align*}
The unitary matrix $\mI - \tau\vv\vv^*$ is called a ``Householder reflection'' \cite[sec.\ 5.1.2]{golub_van_loan}.

Because $\mQ$ is implicitly represented through $\vv_1,\, \ldots,\, \vv_k$ and $\tau_1,\, \ldots,\, \tau_k$, left-multiplying a matrix $\mX$ by $\mQ$ or $\mQ^*$ involves $k$ updates of the form $\mX \leftarrow \mX - \tau_i\vv_i(\vv_i^*\mX)$, with $\vv_i^*\mX$ computed by level-2 BLAS operations. Vectorization and multithreading considerations mean that a single $(k \times m) \cdot (m \times n)$ product with level-3 BLAS is significantly faster than $k$ sequential $(1 \times m) \cdot (m \times n)$ products with level-2 BLAS \cite[sec.\ 1.5]{golub_van_loan}. For this reason, it is advantageous to define $\mV = [\vv_1 \: \cdots \: \vv_k] \in \C^{m \times k}$ and to seek a $\mT \in \C^{k \times k}$ such that
\begin{equation}
\mQ = \mI - \mV\mT\mV^*, \label{eqn:compact_wy}
\end{equation}
so that $\mQ$ can be efficiently applied using $\mQ\mX = \mX - \mV\mT(\mV^*\mX)$, with $\mV^*\mX$ evaluated using level-3 BLAS routines.

If $\mV$ is lower-triangular with unit diagonal and $\mT$ is upper-triangular, then we say that \cref{eqn:compact_wy} constitutes a \emph{compact WY form} of $\mQ$ \cite{schreiber_van_loan_compactwy}. Representing $\mQ$ in compact WY form is straightforward when $(\vv_1,\, \tau_1),\, \ldots,\, (\vv_k,\, \tau_k)$ are generated by \cref{alg:golub_businger}, for two reasons. First, the upper-triangular structure of $\mR$ already ensures that $\mV = [\vv_1 \:\cdots\: \vv_k]$ is lower-triangular with unit diagonal. Second, regardless of the structure $\mV$, an upper-triangular $\mT$ satisfying \cref{eqn:compact_wy} always exists and can be easily generated using an algorithm developed by Schreiber and Van Loan \cite{schreiber_van_loan_compactwy}.

\subsection{Accelerated CPQR Algorithms} \label{subsec:blocked_cpqr_review}
When the dimensions of $\mA$ are extremely large, \cref{alg:golub_businger} suffers a computational bottleneck that results from sequential applications of Householder reflections using level-2 BLAS operations. For non-pivoted QR factorizations, it is standard practice to limit level-2 BLAS work by manipulating columns in blocks, where each block yields a set of Householder vectors that are accumulated into compact WY form for later application to the rest of the matrix using level-3 BLAS \cite[sec 5.2.3]{golub_van_loan}. In \cref{alg:golub_businger}, however, Householder reflectors are determined on the basis of a column swap which requires all previous reflections to be applied to the matrix \emph{before} a compact WY representation can be formed. In this setting, limiting the amount of work performed at level-2 BLAS is more difficult.

This paper describes a modification of \cref{alg:golub_businger} that addresses this bottleneck for matrices with certain structure, but we are by no means the first authors to tackle this problem. A large class of algorithms have been proposed that reduce the number of rows manipulated by level-2 BLAS, such as the following.

\begin{itemize}[leftmargin=*,labelsep=0pt]
\item[] \textbf{Partial Householder reflections.} Selecting pivot columns in line \eqref{line:GEQP3_argmax} of \cref{alg:golub_businger} requires knowledge of column norms in the $(i\mcol m,i \mcol n)$ block of $(\mI - \tau_i\vv_i\vv_i^*)\mR$ as computed in line \eqref{line:GEQP3_R_update}, but it is well known that only the first row of this matrix is needed to determine these norms \cite[sec 5.4.2]{golub_van_loan}. Quianta-Ort\'i, Sun, and Bishof \cite{orti_blas3_cpqr} use this fact to limit level-2 BLAS work to a small number of rows. Their method is to select pivot columns in blocks where, at each step, a Householder reflection is applied ``partially'' to reveal the single row needed to update column norms. Meanwhile, reflections are accumulated into a compact WY form that allows the remaining rows to be updated with level-3 BLAS at the end of the block. This strategy is used in GEQP3, the LAPACK implementation of \cref{alg:golub_businger}.
\item[] \textbf{Sketching.} Many tasks in numerical linear algebra can be accelerated by first generating a random ``sketching matrix'' $\mOmega \in \R^{m \times d}$, with $d \ll m$ \cite{halko_finding_structure_with_randomness}. If drawn from a suitable distribution, $\mOmega$ induces a transform $\vx \mapsto \mOmega^*\vx$ that randomly embeds vectors in a much lower dimensional space, while approximately preserving geometry in a sense made precise by the Johnson-Lindenstrauss lemma \cite[Thm 1.6.1]{vershynin_prob_methods_for_data}. Duersch and Gu \cite{duersch_gu_rqrcp} used this fact devise an ``RCPQR'' algorithm which selects pivot columns using a CPQR factorization of $\mOmega^*\mA$, where the entries of $\mOmega$ are i.i.d.\ standard Gaussians. This limits BLAS-2 computation to the much smaller number of rows in $\mOmega^*\mA$, and once pivot columns are selected, the remaining transformations necessary to CPQR-factorize $\mA$ proceed with BLAS-3. A related class of ``randomized interpolative decomposition'' (RID) algorithms use sketching with non-Gaussian $\mOmega$ to form a column-based low-rank approximation of $\mA$ \cite{martinsson_rokhlin_tygert_randid, woolfe_srft}.

\item[] \textbf{Randomized blocking.} Because RCPQR selects at most $r \defeq \rank \mOmega^*\mA$ skeleton columns, the constraint $r \leq \min\{ d,\, n \}$ means that modifications are needed when a full set of $\min\{ m,\, n \}$ skeleton columns is desired. In this vein, Martinsson et al.\ \cite{martinsson_randomized_qrcp} describe a variant of \cref{alg:golub_businger} which selects a full skeleton set in blocks of $b$, where each block is selected by applying RCPQR with $d = b$ to the trailing submatrix of $\mR$. The Householder reflections generated for each block can be accumulated into compact WY form and applied to the rest of the matrix with BLAS-3. A similar blocking strategy is considered in \cite{duersch_gu_rqrcp}.

\item[] \textbf{GKS and DEIM.} If only $k \ll n$ skeleton columns are needed, the GKS \cite{golub_klema_stewart} and DEIM \cite{sorensen_embree_deim} algorithms provide a means to select these columns by pivoting on a matrix of only $k$ rows, corresponding to the leading $k$ right singular vectors of $\mA$. Computing these vectors is, in and of itself, a computationally intensive task. For this reason, randomized GKS \cite{rgks} or randomized DEIM \cite{saibaba_rdeim} algorithms that estimate the leading singular vectors with sketching may be preferable.
\end{itemize}

This paper is interested in matrices that have far more columns than rows, such that the techniques above do not address the main bottleneck. Matrices of this aspect ratio are amenable to a separate class of algorithms that leverage parallelization or randomization. These include the following.
\begin{itemize}[leftmargin=*,labelsep=0pt]
\item[] \textbf{Global and local pivoting.} Bishof \cite{bishof_local_pivoting} describes two modifications of \cref{alg:golub_businger} that parallelize the column load across several processors, a ``global'' version and a ``local'' one. The ``global'' algorithm distributes column norm updates across different processors, which each send updated norms to a ``lead'' processor with authority to choose the next pivot. To avoid the bottleneck that results from each processor waiting for the lead's decision, the ``local'' pivoting method equips each processor with an incremental condition number estimator that allows it to select pivot columns on its own. Both strategies employ row-restricted Householder reflections, as in \cite{orti_blas3_cpqr}, to allow norm updates and pivot selections to be performed more-or-less simultaneously.

\item[] \textbf{Tournament pivoting.} Demmel et al.\ \cite{demmel_tournament_pivoting} developed a parallelized CPQR factorization algorithm that selects pivots in blocks of $b$, while minimizing the communication cost associated with column movement for each block. This is accomplished in a multi-level ``tournament,'' where at the first level, each processor is assigned a small block of columns from which it selects $b$ skeleton columns using a rank-revealing QR (RRQR) factorization. At the next level, each processor is assigned a combined block of $2b$ skeleton columns which is again downsampled to $b$ columns by an RRQR factorization. After $\cO(\log n)$ rounds of the tournament, a block of $b$ columns has been selected; this is repeated until a full skeleton set is obtained.

\item[] \textbf{Random sampling.} While they are not CPQR algorithms \emph{per se}, there is a large class of algorithms developed in theoretical computer science that choose column subsets by randomly sampling from an appropriate distribution. The first algorithms in this category sampled based on column norm alone \cite{frieze_monte_carlo}, while later algorithms sampled based on subspace leverage leverage scores \cite{mahoney_drineas_muthukrishnan, mahoney_cur}. These methods require significant oversampling. A related class of ``volume sampling'' algorithms sample columns in batches that bias towards high linear independence \cite{deshpande_rademacher_volumesampling, deshpande_rademacher_projectiveclustering}. These do not require the same degree of oversampling, but the samples themselves are more expensive to generate. Unique among this class of algorithms is the SE-QRSC algorithm of Fakih and Grigori \cite{fakih_grigori_sparse_embedding}. This algorithm applies CPQR-based column selection to $\mA\mPhi$, where $\mPhi$ is a ``sparse embedding'' matrix that linearly combines small groups of columns from $\mA$. A CPQR factorization of $\mA\mPhi$ yields a set of ``skeleton column groups,'' which can then be mapped back to skeleton columns of the original matrix.
\end{itemize}

These algorithms are effective at rapidly selecting columns from matrices with far more columns than rows. In \cref{section:cceqr}, we will focus on accelerating the Golub-Businger algorithm (\cref{alg:golub_businger}) for matrices of this sort. \cref{alg:golub_businger} is distinguished from the techniques just described by the fact that its column choice satisfies a readily interpretable greedy ordering, namely, that of $\mathrm{GB}(k)$ form. In practice, this does not make \cref{alg:golub_businger} perform significantly better or worse than the other algorithms described here in the sense of optimizing the rank-revealing bound in \cref{eqn:rrqr_factorization_1,eqn:rrqr_factorization_2}.
\section{The CCEQR Algorithm} \label{section:cceqr}
\begin{figure}
    \centering
    \includegraphics[scale=.8]{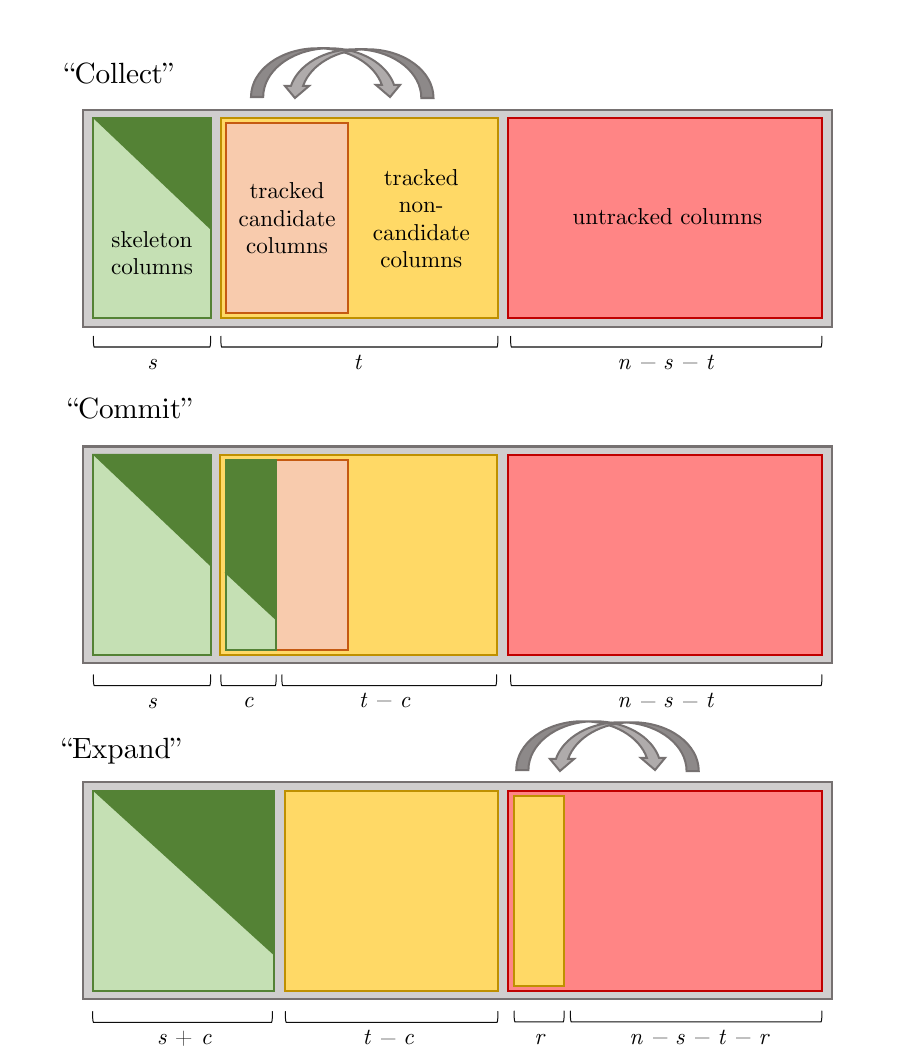}
    \caption{A schematic representation of CCEQR. In the ``collect'' stage, a set of candidate skeleton columns are selected from the tracked set. The ``commit'' stage chooses a subset of these candidates to bring into the skeleton. The ``expand'' stage brings new columns into the tracked set.} \label{fig:cceqr_schematic}
\end{figure}

CCEQR works in cycles, each of which pivots at least one column into the skeleton set. Each cycle starts with a permutation $\mPi$ and a unitary $\mQ$ such that $\mQ^*\mA\mPi$ has $\mathrm{GB}(s)$ form for some $s \geq 0$. The columns of $\mA\mPi$ are partitioned into three contiguous blocks:
\begin{enumerate}
\item The leftmost block contains $s$ ``committed'' skeleton columns. These will not be modified by the present cycle or by any future cycles.
\item The middle block contains ``tracked'' columns whose norms in the subspace orthogonal to $\mA\mPi(:,1 \mcol s)$ are monitored by the algorithm.
\item The rightmost block contains ``untracked'' columns, for which the only information known to CCEQR is the starting norm.
\end{enumerate}
The objective of a given cycle is to modify $\mPi$ and $\mQ$ such that $\mQ^*\mA\mPi$ has $\mathrm{GB}(s+c)$ form for some $c \geq 1$. This is accomplished by a procedure which can be informally described as follows.
\begin{enumerate}
\item A ``collect'' step assembles a small set of candidate skeleton columns from the tracked set, and forms a CPQR factorization of the candidates.
\item A ``commit'' step uses the CPQR factors to identify $c$ columns that can be safely added to the skeleton, and modifies $\mPi$ and $\mQ$ accordingly. If this action brings the skeleton to size $k$, then the algorithm terminates.
\item Otherwise, an ``expand'' step moves a portion of the untracked columns into the tracked set, adding new candidates to be collected at the next cycle.
\end{enumerate}
After at most $k$ cycles, $\mQ^*\mA\mPi$ will acquire GB$(k)$ form. This cycle is represented schemtically in \cref{fig:cceqr_schematic}.

Two major factors make this an efficient procedure for column selection. First, BLAS-2 Householder reflections are only used during the CPQR factorization in the ``collect'' step, during which they are only applied to small minority of columns (i.e., the candidate set). Second, and more significantly, Householder reflections in the ``commit'' and ``expand'' steps (performed with BLAS-3) are applied only to the set of tracked columns. This can greatly improve runtime, because if the parameters of CCEQR are properly tuned (cf.\ \cref{section:experiments}), then the tracked set often represents a very small proportion of the entire matrix.

\subsection{Preliminaries}
At the outset of each cycle, CCEQR has computed a unitary $\mQ \in \C^{m \times m}$ and a permutation $\mPi \in \{ 0,\, 1 \}^{n \times n}$ such that $\mQ^*\mA\mPi$ has GB$(s)$ form for some $s \geq 0$, starting with $\mQ = \mI_{m \times m},\, \mPi = \mI_{n \times n}$, and $s = 0$ at the first cycle. CCEQR applies column permutations globally, but for efficiency, rotations are applied only to the tracked set. Thus, CCEQR does not have explicitly compute $\mQ^*\mA\mPi$, but instead computes
\begin{equation*}
\mR \defeq \begin{bNiceMatrix}[last-row]
\mR_1 & \mR_2 & \mX \\
\mbox{\scriptsize $s$} & \mbox{\scriptsize $t$} & \mbox{\scriptsize $n-s-t$}
\end{bNiceMatrix}
\end{equation*}
where $\mR_1$ is an upper-triangular matrix corresponding to $s$ ``committed'' skeleton columns, $\mR_2$ contains $t$ ``tracked'' columns, and the remaining ``untracked'' columns make up $\mX$. We have
\begin{equation*}
\begin{bmatrix} \mR_1 & \mR_2 \end{bmatrix} = \mQ^*\mA\mPi(:,1 \mcol (s+t)) \quad\text{and}\quad \mX = \mA\mPi(:,(s+t+1) \mcol n).
\end{equation*}
CCEQR maintains a vector $\vgamma$ that records the \emph{residual} norms of tracked columns, and the \emph{overall} norms of untracked columns:
\begin{equation*}
\vgamma(j) = \begin{cases}
\| \proj_{(\range \mA\mPi(:,1 \mcol s))^\perp} \mA\mPi(:,j) \|_2^2 & j \leq s+t \\
\| \mA\mPi(:,j) \|_2^2 & j > s+t
\end{cases}.
\end{equation*}
The partition into committed, tracked, and untracked columns is maintained such that at the beginning of each cycle,
\begin{equation}
\max_{s < j \leq s+t} \vgamma(j) \geq \max_{s + t < j} \vgamma(j), \label{eqn:cceqr_norm_bounds}
\end{equation}
a relation which will ensure that each cycle of CCEQR adds at least one new skeleton column. Equation \cref{eqn:cceqr_norm_bounds} holds automatically for the first cycle, which is initialized with $s = 0$ and $t = n$ (following the convention that $|\mR(0,0)| = \infty$).

CCEQR also maintains a lower-triangular $\mV$ and upper-triangular $\mT$ providing a compact WY representation of $\mQ$. These matrices are padded with zeros so that $\mV \in \C^{m \times k}$ and $\mT \in \C^{k \times k}$, and the first cycle initializes with $\mV = \mZero_{m \times k},\, \mT = \mZero_{k \times k}$. Lastly, $\mPi$ is represented via a permutation vector $\vp \in [n]^n$ such that $\mA\mPi = \mA(:,\vp)$; this vector is modified by each cycle.

\Cref{alg:cceqr_initialize} summarizes, in pseudocode, the steps taken to initialize the first cycle of CCEQR. Note that $\mu$ stores the value of $\max_{s + t<j} \vgamma(j)$.
\begin{algorithm}
\caption{\textsc{CCEQR: Subroutine ``Initialize''}} \label{alg:cceqr_initialize}
\begin{algorithmic}[1]
\STATE \textbf{inputs:} $\mA,\, k$.
\item[]
\STATE $\vp \leftarrow [1,\, 2,\, \ldots,\, n],\, \mR \leftarrow \mA,\, \mV \leftarrow \mZero_{n \times k}$, and $\mT \leftarrow \mZero_{k \times k}$.
\STATE $s \leftarrow 0,\, t \leftarrow n,\, \mu \leftarrow 0$.
\item[]
\FOR{$j = 1,\, \ldots,\, n$}
\STATE $\vgamma(j) \leftarrow \| \mR(:,j) \|_2^2$.
\ENDFOR
\item[]
\RETURN $s,\, t,\, \mu,\, \vgamma,\, \vp,\, \mV,\, \mT,\, \mR$.
\end{algorithmic}
\end{algorithm}

\subsection{``Collect'' Step} \label{subsec:collect}
This begins by forming a block of $b$ ``candidate'' skeleton columns from the tracked set, where $b = 1 + \lfloor \rho (t-1) \rfloor$ and $\rho \in (0,1)$ is a fixed parameter set by the user\footnote{We define $b$ in this way so that $1 \leq b < t$ whenever $t > 1$, and otherwise $b = t = 1$.}. Typically we use $\rho \ll 1$ so that $b \ll t$. For the first cycle of CCEQR, which initializes with $t = n$, we follow this by setting $t \leftarrow b + 1$.

The candidate block consists of the tracked columns with largest residual norm, i.e., column indices
\begin{equation*}
s + \vsigma(1 \mcol b) \defeq [s + \vsigma(1),\, \ldots,\, s + \vsigma(b)]
\end{equation*}
where $\vsigma$ is a permutation vector that sorts $\vgamma((s+1) \mcol (s+t))$ in descending order. After forming the candidate block, we must gather information that will be necessary for deciding which columns can be safely added to the skeleton in the ``commit'' step. To that end we record\footnote{If all indices in $[s+1,\, s+t]$ are candidates, which occurs only if $b = t = 1$, then $\delta \defeq 0$.}
\begin{equation*}
\delta \defeq \max\{ \vgamma(j) \suchthat j \text{ is a non-candidate index in } [s+1,\, s+t] \},
\end{equation*}
and we form a CPQR factorization of the candidate columns' residuals:
\begin{align}
&\mB \leftarrow \mR((s+1) \mcol m,\, s+\vsigma(1 \mcol b))  \nonumber \\
&\widehat{\vp},\, \widehat{\vtau},\, \widehat{\mV},\, \widehat{\mR} \leftarrow \textsf{GEQP3}(\mB), \label{eqn:cceqr_geqp3}
\end{align}
where $\vb$ is a permutation of $1 \mcol b$, $\widehat{\vtau} \in \R^d$ and $\widehat{\mV} \in \C^{(n - s) \times d}$ are Householder scalars and vectors (with $d = \min\{ m-s,\, b \}$), and $\widehat{\mR}$ is the resulting R-factor. To finish the ``collect'' step, we permute the candidate columns to indices $(s+1) \mcol (s+b)$ of $\mR$ and reorder them according to $\widehat{\vp}$. We also permute $\vp$ and $\vgamma$ for consistency. Because all column permutations up to this point have taken place within tracked set, equation \cref{eqn:cceqr_norm_bounds} is still in force when the ``collect'' step terminates.

\Cref{alg:cceqr_collect} summarizes the ``collect'' step in pseudocode, where
\begin{equation*}
\widehat{\mX} \leftarrow \textsf{PermuteColumns}(\mX,\, \vy,\, \vz,\, \vw)
\end{equation*}
reorders the columns of $\mX$ such that $\widehat{\mX}(:,\vz) = \mX(:,\vy(\vw))$. $\textsf{PermuteEntries}$ reorders the entries of a vector in an identical fashion.
\begin{algorithm}
\caption{\textsc{CCEQR: Subroutine ``Collect''}} \label{alg:cceqr_collect}
\begin{algorithmic}[1]
\STATE \textbf{inputs:} $\rho,\, s,\, t,\, \vgamma,\, \vp,\, \mR$.
\item[]
\STATE $b \leftarrow 1 + \lfloor \rho(t - 1) \rfloor$.
\STATE $\vsigma \leftarrow$ permutation of $1 \mcol t$ that sorts $\vgamma((s+1) \mcol (s+t))$ in decreasing order.
\STATE $\delta \leftarrow \vgamma(s + \vsigma(b + 1))$ \textbf{if} $b < t$, \textbf{else} 0.
\STATE \textbf{if }\textsc{first-cycle}\textbf{ then } $t \leftarrow b$.
\item[]
\STATE $\widehat{\vp},\, \widehat{\vtau},\, \widehat{\mV},\, \widehat{\mR} \leftarrow \textsf{GEQP3}(\mR((s+1) \mcol m,\, s + \vsigma(1 \mcol b)))$. \label{line:collect_geqp3}
\STATE $\mR \leftarrow \textsf{PermuteColumns}(\mR,\, s+\vsigma(1 \mcol b),\, (s+1) \mcol (s+b),\, \widehat{\vp})$. \label{line:collect_permutation}
\STATE $\vgamma \leftarrow \textsf{PermuteEntries}(\vgamma,\, s+\vsigma(1 \mcol b),\, (s+1) \mcol (s+b),\, \widehat{\vp})$.
\STATE $\vp \leftarrow \textsf{PermuteEntries}(\vp,\, s+\vsigma(1 \mcol b),\, (s+1) \mcol (s+b),\, \widehat{\vp})$. \label{line:cceqr_collect_finalperm}
\item[]
\RETURN $\delta,\, \widehat{\vtau},\, \widehat{\mV},\, \widehat{\mR}$.
\end{algorithmic}
\end{algorithm}

\subsection{``Commit'' Step} \label{subsec:commit}
The goal in this step is to modify $\mQ$ and $\mPi$ based on information prepared during the ``collect'' step such that $\mQ^*\mA\mPi$ has $\mathrm{GB}(s+c)$ form for some $c \geq 1$. Following this, columns $1 \mcol (s+c)$ will count as ``committed'' skeleton columns, in the sense that no future cycle will modify or reorder them. The appropriate changes to $\mQ$ and $\mPi$ are determined using \cref{lemma:commit_rule}.
\begin{lemma} \label{lemma:commit_rule}
Let $\delta,\, \widehat{\vtau},\, \widehat{\mV},\, \widehat{\mR}$ be the data returned by the ``collect'' step of CCEQR (\cref{alg:cceqr_collect}), and let $\widehat{\mPi}$ be the column permutation matrix applied in line \eqref{line:collect_permutation} of that step. Let $\mu = \max_{s+t < j} \| \mA\mPi(:,j) \|_2^2$, and define
\begin{equation}
c = \max\{ i \geq 1 \suchthat |\widehat{\mR}(i,i)|^2 \geq \delta \lor \mu \}. \label{eqn:cceqr_acceptance_number}
\end{equation}
Let $\widehat{\mQ}$ be the unitary matrix that applies the first $c$ Householder reflections computed in the ``collect'' step to the bottom $n - s$ rows, i.e.,
\begin{equation*}
\widehat{\mQ} = \begin{bmatrix}
\mI_{s \times s} & \mZero_{s \times (n-s)} \\
\mZero_{(n-s) \times s} & \prod_{i = 1}^c (\mI - \vtau(i)\widehat{\mV}(:,i)\widehat{\mV}(:,i)^*)
\end{bmatrix}.
\end{equation*}
Then, $(\mQ\widehat{\mQ})^*\mA(\mPi\widehat{\mPi})$ has $\mathrm{GB}(s+c)$ form.
\end{lemma}

\begin{proof}
A detailed proof is given in \cref{section:commit_rule_proof}. A proof sketch is as follows: \cref{eqn:cceqr_norm_bounds} guarantees that the maximum in \cref{eqn:cceqr_acceptance_number} is over a nonempty set, so $c \geq 1$ is well-defined. To bring $\mR$ to $\mathrm{GB}(k)$ form in the usual way, we would apply \cref{alg:golub_businger} to $\mB_1 \defeq \mR((s+1) \mcol m,\, (s+1) \mcol n)$. In the ``collect'' step, we have instead applied \cref{alg:golub_businger} to $\mB_2 \defeq \mR((s+1) \mcol m,\, (s+1) \mcol (s+b))$. The resulting column permutation brings $\mR$ to $\mathrm{GB}(s+c_0)$ form, where $c_0$ is the number of iterations before \cref{alg:golub_businger} applied to $\mB_1$ selects a different column than it would if applied to $\mB_2$. This occurs when there is a non-candidate column whose residual norm (if included in pivoting) would be bigger than that of any of the candidates. Equation \cref{eqn:cceqr_acceptance_number} causes CCEQR to only ``commit'' a pivot columns if its residual is greater than the \emph{overall} norm of all non-candidates. Thus ensuring no unwanted columns are accepted.
\end{proof}

As noted in the statement of \cref{lemma:commit_rule}, the permutation $\mPi \mapsto \mPi\widehat{\mPi}$ was already performed by the ``collect'' step at line \eqref{line:collect_permutation}. It remains to update $\mQ \mapsto \mQ\widehat{\mQ}$ by incorporating the first $c$ Householder reflectors from the ``collect'' step into the compact WY form for $\mQ$. To that end, the Schreiber-van Loan algorithm \cite{schreiber_van_loan_compactwy} is used to form an upper-triangular $\widehat{\mT} \in \C^{c \times c}$ such that
\begin{equation*}
\prod_{i = 1}^c (\mI - \vtau(i)\widehat{\mV}(:,i)\widehat{\mV}(:,i)^*) = \mI - \widehat{\mV}(:,1 \mcol c)\widehat{\mT}\widehat{\mV}(:,1 \mcol c)^*.
\end{equation*}
In prepration for residual column norm updates, this rotation is applied to the bottom $n-s$ rows of the tracked set. We then update $\mQ$ by incorporating $\widehat{\mT}$ and $\widehat{\mV}(:,1 \mcol c)$ into the global compact WY factors $\mT$ and $\mV$. Although we could incorporate each reflector serially using the Schreiber-Van Loan algorithm \cite{schreiber_van_loan_compactwy}, we instead apply a block update to the WY factors using a small number of BLAS-3 operations. The details of this are explained in \cref{section:efficient_wy_details}.

Residual column norms are now updated in the tracked set by computing
\begin{equation}
\vgamma(j) \leftarrow \vgamma(j) - \| \mR((s+1) \mcol (s+c),\, j) \|_2^2 \label{eqn:cceqr_colnorm_update}
\end{equation}
for $j = s+c+1,\, \ldots,\, s+t$. Note that $\| \mR((s+1) \mcol (s+c),\, j) \|_2^2$ represents the portion of $\mR(:,j)$ spanned by the newly committed skeleton columns. We record these new commits by setting $s \leftarrow s+c$ and $t \leftarrow t-c$. The numerical stability of the column-norm update in \cref{eqn:cceqr_colnorm_update} could be improved using techniques from \cite{zlatko_bujanovic_rrqr_failure}, though we do not do this here.

\Cref{alg:cceqr_commit} summarizes the ``commit'' step. In this pseudocode, the function $\mX \leftarrow \textsf{CompactWY}(\vx,\, \mY)$ takes as input Householder scalars $\vx$ and reflectors $\mY$, and returns the upper-triangular $\mX$ needed to represent the corresponding unitary matrix in compact WY form. The function $\widehat{\mS} \leftarrow \textsf{ApplyQt}(\mS,\, \mX,\, \mY,\, \vi,\, \vj)$ applies a compact WY-form rotation $(\mI - \mY\mX\mY^*)^*$ to the $(\vi,\, \vj)$ block of $\mS$ using BLAS-3 routines\footnote{To fully take advantage of BLAS-3 alongside the triangular structures inherent to compact WY form, our implementation of this function is templated off of the LAPACK routine LARFB.}. Lastly, $\textsf{UpdateWY}$ adds new Householder reflectors to $\mQ$ using formulas in \cref{section:efficient_wy_details}.
\begin{algorithm}
\caption{\textsc{CCEQR: Subroutine ``Commit''}} \label{alg:cceqr_commit}
\begin{algorithmic}[1]
\STATE \textbf{inputs:} $\delta,\, \mu,\, s,\, t,\, \vgamma,\, \widehat{\vtau},\, \widehat{\mV},\, \mT,\, \mV,\, \widehat{\mR},\, \mR$.
\item[]
\STATE $c \leftarrow \max\{ i \geq 1 \suchthat |\widehat{\mR}(i,i)|^2 \geq \delta \lor \mu \}$. \label{line:cceqr_acceptance_size}
\STATE $\widehat{\mT} \leftarrow \textsf{CompactWY}(\widehat{\vtau}(1:c),\, \widehat{\mV}(:,1 \mcol c))$.
\STATE $\mR \leftarrow \textsf{ApplyQt}(\mR,\, \widehat{\mT},\, \widehat{\mV}(:,1:c),\, (s+1) \mcol m,\, (s+1) \mcol (s+t))$.
\STATE $\mT,\, \mV \leftarrow \textsf{UpdateWY}(\mT,\, \mV,\, \widehat{\mT},\, \widehat{\mV}(:,1 \mcol c))$.
\item[]
\item $M \leftarrow 0$.
\FOR{$j = (s+1) \mcol (s+t)$}
\STATE $\vgamma(j) \leftarrow \vgamma(j) - \| \mR((s+1) \mcol (s+c),\, j) \|_2^2$.
\STATE $M \leftarrow \max(M,\, \vgamma(j))$. \label{line:cceqr_commit_maxtrackednorm}
\ENDFOR
\item[]
\STATE $s \leftarrow s+c$.
\STATE $c \leftarrow c-t$.
\item[]
\RETURN $M$.
\end{algorithmic}
\end{algorithm}

\subsection{``Expand'' Step} \label{subsec:expand}
This step executes if the ``commit'' step did not complete the skeleton (i.e., $s+c < k$) and there are still untracked columns (i.e., $s+t < n$). In this case, because the tracked columns were orthogonalized against new skeleton columns, it is possible that their residual norms have decreased leading to a violation of inequality \cref{eqn:cceqr_norm_bounds}. To ensure that the next cycle is able to commit at least $1$ additional skeleton column, we must modify the partition into ``tracked'' and ``untracked'' columns such that $\cref{eqn:cceqr_norm_bounds}$ is restored.

To that end, the ``expand'' step assembles an index vector $\vu$ containing all untracked indices $j$ such that $\| \mA\mPi(:,j) \|_2^2 \geq \max_{s< j \leq s+t} \vgamma(j)$.These indices will be moved from untracked to tracked\footnote{In the rare event that no untracked columns meet the norm threshold, we lower it to $0.9\max_{s < j \leq s+t}$ to ensure that the tracked set expands.}. While assembling $\vu$,\, we also record
\begin{equation*}
\mu \leftarrow \max\{ \| \mA\mPi(:,j) \|_2^2 \suchthat j > s+t,\, \| \mA\mPi(:,j) \|_2^2 < \max_{s< j \leq s+t} \vgamma(j) \},
\end{equation*}
which will represent the maximum untracked column norm at the next cycle\footnote{This is needed in line \eqref{line:cceqr_acceptance_size} of the ``commit'' step (\cref{alg:cceqr_commit}); see also \cref{lemma:commit_rule}.}. The columns associated with $\vu$ are permuted to positions $(s+t+1) \mcol (s+t+r)$ (where $r = \textsf{length}(\vu)$) and orthogonalized against the skeleton. We record their new residual norms in $\vgamma$ and we set $t \leftarrow t + r$, thus completing the expansion of the ``tracked'' set.

\cref{alg:cceqr_expand} summarizes the ``expand'' step, where $\vj,\, \beta \leftarrow \textsf{Threshold}(\vx,\, i,\, \alpha)$ returns a vector $\vj$ containing all indices $j \geq i$ such that $\vx(j) \geq \alpha$, along with the value $\beta = \max_{j \geq i,\, j \not\in \vj} \vx(j)$. For the other functions used in this pseudocode, refer to \cref{subsec:collect,subsec:commit}.
\begin{algorithm}
\caption{\textsc{CCEQR: Subroutine ``Expand''}} \label{alg:cceqr_expand}
\begin{algorithmic}[1]
\STATE \textbf{inputs:} $s,\, t,\, M,\, \vgamma,\, \mT,\, \mV,\, \mR$.
\item[]
\STATE $\vu,\, \mu \leftarrow \textsf{Threshold}(\vgamma,\, s+t+1, M)$.
\STATE $r \leftarrow \mathsf{length}(\vu)$.
\IF{$r = 0$}
\STATE $\vu,\, \mu \leftarrow \textsf{Threshold}(\vgamma,\, s+t+1, 0.9\mu)$
\STATE $r \leftarrow \mathsf{length}(\vu)$.
\ENDIF
\item[]
\STATE $\mR \leftarrow \mathsf{PermuteColumns}(\mR,\, \vu,\, (s+t+1) \mcol (s+t+r),\, 1 \mcol r)$.
\STATE $\vgamma \leftarrow \mathsf{PermuteEntries}(\vgamma,\, \vu,\, (s+t+1) \mcol (s+t+r),\, 1 \mcol r)$.
\STATE $\vp \leftarrow \mathsf{PermuteEntries}(\vp,\, \vu,\, (s+t+1) \mcol (s+t+r),\, 1 \mcol r)$.
\STATE $\mR \leftarrow \mathsf{ApplyQt}(\mR,\, \mT,\, \mV,\, 1 \mcol m,\, (s+t+1) \mcol (s+t+r))$.
\item[]
\FOR{$j = (s+t+1) \mcol (s+t+r)$}
\STATE $\vgamma(j) \leftarrow \vgamma(j) - \| \mR(1 \mcol s,\, j) \|_2^2$.
\ENDFOR
\item[]
\STATE $t \leftarrow t + r$.
\end{algorithmic}
\end{algorithm}

\subsection{Full Algorithm}
For completeness, we provide full pseudocode for CCEQR in \cref{alg:cceqr}. Note that \cref{alg:cceqr} allows for a ``CSSP only'' mode of CCEQR, where only selected columns of $\mR$ are formed, as well as a ``full CPQR'' version where all Householder reflections are applied to the untracked set to produce all of $\mR$. This is controlled by the ``if-then'' block starting at line \eqref{line:cceqr_full_switch}.
\begin{algorithm}
\caption{\textsc{Collect-Commit-Expand QR (CCEQR)}} \label{alg:cceqr}
\begin{algorithmic}[1]
\STATE \textbf{inputs:} $\mA \in \C^{m \times n},\, k \leq \min\{ m,\, n \},\, \rho \in (0,\, 1)$, and $\textsc{full} \in \{ \textsc{True},\, \textsc{False} \}$.
\item[]
\STATE $s,\, t,\, \mu,\, \vgamma,\, \vp,\, \mV,\, \mT,\, \mR \leftarrow \textsf{Initialize}(\mA,\, k)$.\qquad\emph{\# see \cref{alg:cceqr_initialize}}
\item[]
\WHILE{$s < k$}
\STATE $\delta,\, \widehat{\vtau},\, \widehat{\mV},\, \widehat{\mR} \leftarrow \textsf{Collect}(\rho,\, s,\, t,\, \vgamma,\, \vp,\, \mR)$.\qquad\emph{\# see \cref{alg:cceqr_collect}}
\STATE $M \leftarrow \textsf{Commit}(\delta,\, \mu,\, s,\, t,\, \vgamma,\, \widehat{\vtau},\, \widehat{\mV},\, \mT,\, \mV,\, \widehat{\mR},\, \mR)$.\qquad\emph{\# see \cref{alg:cceqr_commit}}
\item[]
\IF{$s < k$}
\STATE $\textsf{Exand}(s,\, t,\, M,\, \vgamma,\, \mT,\, \mV,\, \mR)$.\qquad\emph{\# see \cref{alg:cceqr_expand}}
\ELSE
\STATE \textbf{goto} line \eqref{line:cceqr_return}.
\ENDIF
\ENDWHILE
\item[]
\IF{$\textsc{full} = \textsc{True}$} \label{line:cceqr_full_switch}
\STATE $\textsf{ApplyQt}(\mR,\, \mT,\, \mV,\, 1 \mcol m,\, (s+t+1) \mcol n)$.
\ENDIF
\item[]
\RETURN $\vp$. \label{line:cceqr_return}
\end{algorithmic}
\end{algorithm}

We have claimed that CCEQR computes an equivalent column permutation to the Golub-Businger algorithm. We end this section with a formal statement and proof of this claim.
\begin{theorem} \label{theorem:cceqr_equivalence}
Let $\vp \in [n]^k$ be the permutation vector returned by \cref{alg:cceqr}. This permutation is equivalent to the one computed by \cref{alg:golub_businger}, in the sense that there is a unitary $\mQ$ such that $\mQ^*\mA(\mcol,\vp)$ is in $\mathrm{GB}(k)$ form. 
\end{theorem}

\begin{proof}
The proof is inductive. Let $\mQ^{(t)}$ and $\mPi^{(t)}$ be the global column permutation and unitary factor at the end of the $t\nth$ cycle of \cref{alg:cceqr}, with $\mQ^{(0)} = \mI_{m \times m}$ and $\mPi^{(0)} = \mI_{n \times n}$. Then $(\mQ^{(0)})^*\mA\mPi^{(0)}$ has $\mathrm{GB}(s_0)$ form with $s_0 = 0$. Suppose that $(\mQ^{(t-1)})^*\mA\mPi^{(t-1)}$ has $\mathrm{GB}(s_{t-1})$ form for some $t \geq 1$. The $t\nth$ cycle transforms $\mQ^{(t-1)} \mapsto \mQ^{(t)}$ and $\mPi^{(t-1)} \mapsto \mPi^{(t)}$ as described in \cref{lemma:commit_rule}, and the conclusion of that lemma is that $(\mQ^{(t)})^*\mA\mPi^{(t)}$ has $\mathrm{GB}(s_{t-1}+c)$ form for some $c \in [1, k - s_{t-1}]$. For some $T \leq k$ we will therefore have $s_T = k$, at which point the algorithm terminates. Hence, the theorem holds with $\mQ = \mQ^{(T)}$ and $\mA(\mcol,\vp) = \mA\mPi^{(T)}$.
\end{proof}
\section{Experiments} \label{section:experiments}
We now present experiments comparing the performance of CCEQR (\cref{alg:cceqr}) and GEQP3 (the implementation of \cref{alg:golub_businger} in LAPACK). We have selected our test cases to illustrate that the comparison between these two algorithms is strongly affected by the structural properties of the matrix being factorized. \Cref{subsec:clustering,subsec:dft} demonstrate examples from scientific applications where CCEQR outperforms GEQP3 by as much as an order of magnitude. \Cref{subsec:random_gaussians} will show that on random unstructured problems, CCEQR and GEQP3 perform more-or-less the same. Finally, \cref{subsec:hadamard_adversary} will demonstrate an ``adversarial'' problem specifically designed to make CCEQR slower than GEQP3.

All experiments in this paper were performed using Julia $1{.}10{.}2$ with the default \texttt{Float64} numerical type. The experiments in this section were performed on a Xeon Platinum 8362 core using the Intel MKL BLAS. Experiments on an Apple M2 chip using both the AppleAccelerate BLAS and OpenBLAS are included in \cref{section:apple_results}.

\subsection{Spectral Demixing} \label{subsec:clustering}
Our first test case involves matrices generated as part of a spectral demixing computation. In this setting we consider data points $\{ x_1,\, \ldots,\, x_n \} \subseteq \cX$ drawn i.i.d.\ from an $m$-component mixture model $\cP$, given by
\begin{equation*}
\cP = \sum_{i = 1}^m \beta_i \cP_i,\quad \beta_1,\, \ldots,\, \beta_m \geq 0,\quad \sum_{i = 1}^m \beta_i = 1,
\end{equation*}
where $\cP_1,\, \ldots,\, \cP_m$ are probability measures. We seek to label each point according to which component it was drawn from using spectral clustering \cite{luxburg_spectral_clustering}. To that end, the CPQR factorization to be computed is
\begin{equation}
\mW_m\tp\mPi = \mQ\mR, \label{eqn:clustering_cpqr}
\end{equation}
where the columns of $\mW_m$ are the leading $m$ eigenvectors of the the normalized kernel evaluation matrix $\mK(i,j) = (d_i d_j)^{-1/2}\cK(x_i,\, x_j)$, with $\cK$ being a positive definite kernel on $\cX$ and $d_i = \sum_{r = 1}^n \cK(x_
i,\, x_r)$. Under suitable conditions on mixture component separation, the leading $m$-dimensional eigenspace of $\mK$ approximates $\cV \defeq \vspan\{ \vv_1,\, \ldots,\, \vv_m \}$, where $\vv_i(j)$ is a kernelized measure of likelihood of $x_j$ under $\cP_i$ \cite{schiebinger_clustering_geometry}. Damle, Minden, and Ying \cite{damle_robust_clustering} demonstrated that points can be labeled by orthogonalizing $\mW_m\tp\mPi(:,1 \mcol m)$ in \cref{eqn:clustering_cpqr} via a polar decomposition, producing a unitary $\mU \in \R^{m \times m}$ that reveals label likelihoods through the relation $\mW_m\mU \approx [\vv_1 \:\cdots\: \vv_m]$.

Our experiments use $\cX = \R^{20}$, and points are drawn from an $m = 20$ component Gaussian mixture model with
\begin{equation*}
\cP_i = \cN(\ell\ve_i,\, \mI),\quad \beta_i = 1,\quad i = 1,\, \ldots,\, m,
\end{equation*}
where $\ve_i$ is the $i\nth$ elementary unit vector and $\ell$ controls the cluster separation scale. We use $\cK(x,\, y) = \exp(-\frac{1}{2\sigma^2}\| x - y \|_2^2)$ for the kernel, with kernel variance $\sigma^2 = 5$. At the scale of this experiment, working directly with $\mK$ is not computationally feasible. We therefore approximate $\mW_m$ by compressing $\mK$ through a pivoted partial Cholesky factorization with significant oversampling, and then directly computing a singular value decomposition of this low-rank approximation. While techniques for handling $\mK$ are not the focus of this paper, matrix-free kernel manipulation methods such as the fast Gauss transform \cite{fast_gauss_transform} would also be appropriate in this setting.

\begin{figure}
    \centering
    \includegraphics[scale=.8]{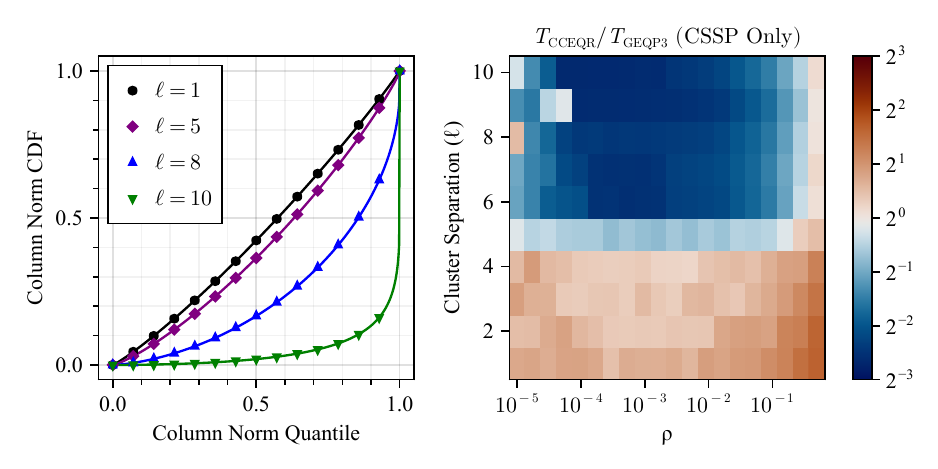}
    \caption{\underline{Left}: cumulative distribution of column norm mass in $\mW_m\tp$ as a function of column norm quantile, versus cluster separation scale $\ell$. \underline{Right}: median runtime ratios for CCEQR and GEQP3 on $m \times n$ matrices generated from spectral demixing with $m = 20$ components and $n = 400{,}000$ data points, across increasing values of $\rho$ and increasing cluster separation, over 30 trials. Cool colors indicate that CCEQR is faster, and warm colors indicate that GEQP3 is faster. Plotted values range from $0{.}16$ to $3{.}24$}. Note that CCEQR was used only to select columns, and the full $\mR$ matrix was not computed. \label{fig:cluster_fixed_size}
\end{figure}
\Cref{fig:cluster_fixed_size} compares the runtimes of CCEQR and GEQP3 over 30 independent trials. For this experiment we varied the cluster separation lengthscale $\ell$, as well as the parameter $\rho$ that controls the size of the candidate column block in CCEQR (cf.\ \cref{subsec:collect}). Increasing $\ell$ means that a larger number of points will be far from the center of any Gaussian component, corresponding to a larger number of small-norm columns in $\mV_m\tp$. This is shown in the left-hand panel of \cref{fig:cluster_fixed_size} which, for several values of $\ell$, plots the percentage of column-norm mass contained at or below a given quantile of the column-norm distribution. Heuristically, we expect CCEQR to run faster than GEQP3 on large problems where the column norm mass is concentrated in a small set of indices (i.e., in the uppermost quantiles). For problems of this sort, the norm-sorting strategy in CCEQR's ``collect'' step is more likely to find all necessary skeleton columns early on, which minimizes the amount of work that must be dedicated to Householder reflections. Indeed, \cref{fig:cluster_fixed_size} shows that increasing $\ell$ leads to a roughly 10x for speedup for CCEQR over GEQP3.

Also apparent in \cref{fig:cluster_fixed_size} are isolated ``islands'' of longer CCEQR runtimes, visible in the upper-left of the right panel. These correspond to a relatively rare scenario wherein, at some cycle in CCEQR, the entire tracked set lies almost entirely in the subspace spanned by the newly committed skeleton columns. In this situation, orthogonalizing against the new skeleton columns will dramatically reduce residual norms in the tracked set, meaning the threshold $M$ computed in the ``commit'' step (\cref{alg:cceqr_commit}) will be quite small. The ``expand'' step (\cref{alg:cceqr_expand}), which adds new tracked columns according to whether or not their norms exceed this threshold, will then make a very large number of columns tracked. The burden of applying Householder reflections to this much larger tracked set greatly slows down future cycles, creating the observed increase in runtimes. In our experience this behavior is more likely to occur when $\rho$ is excessively small, as this produces a small initial tracked set which can be more easily captured in the span of a few skeleton columns.

\begin{figure}
    \centering
    \includegraphics[scale=.8]{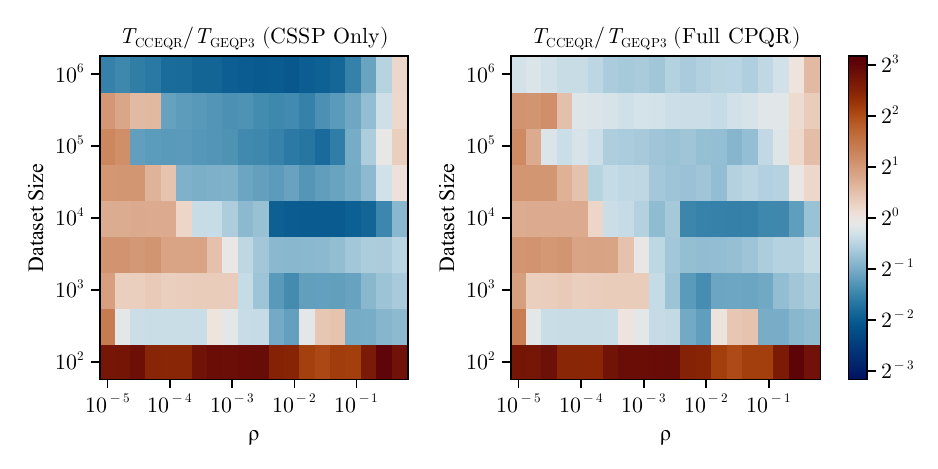}
    \caption{The same experiment as \cref{fig:cluster_fixed_size} with $\ell = 6$ fixed and $n$ increasing, with and without full computation of $\mR$ in CCEQR (note that $\mR$ is not needed for the clustering application). Cool colors indicate CCEQR performing faster than GEQP3, while warm colors indicate GEQP3 performing faster. In the left panel, plotted values range from $0{.}24$ to $8{.}48$. In the right panel, values range from $0{.}34$ to $8{.}56$}. \label{fig:cluster_fixed_scale}
\end{figure}

\Cref{fig:cluster_fixed_scale} shows an identical experiment to \cref{fig:cluster_fixed_size}, except that cluster separation is fixed at $\ell = 6$ and the dataset size is increased over several orders of magnitude. We note that the extra overhead of CCEQR's more complex control flow means it is slower than GEQP3 for problems with relatively few columns. For matrices in this experiment with $10^4$ columns or more, CCEQR is generally faster. This experiment also demonstrates that CCEQR performs best when $\rho$ is neither too large nor too small. Setting $\rho$ too large means CCEQR must deal with excessively large candidate blocks at each cycle. On the other hand, setting $\rho$ to small means that the maximum over tracked residual norms in line \eqref{line:cceqr_commit_maxtrackednorm} of the ``commit'' step (cf.\ \cref{alg:cceqr_commit}) is taken over a much smaller set. Because this maximum is used select new tracked columns, (cf.\ \cref{alg:cceqr_expand}), taking the maximum over too small a set may lead to an unnecessarily large expansion of the tracked set, forcing the algorithm to devote more work to Householder reflections at future cycles. We are pleased to see that for each problem size, the range of ``good'' choices for $\rho$ is fairly broad.

\subsection{Density Functional Theory} \label{subsec:dft}
A fundamental task in computational chemistry is to compute the ground-state electronic energy of a molecular system, which for a system of $m$ elections involves minimizing an energy functional of an $m$-particle wavefunction. Density functional theory (DFT) approximates the $m$-body ground-state wavefunction in terms of $m$ single-particle wavefunctions that correspond to the solution of a nonlinear eigenvalue problem \cite{lin_lu_est}. It is computationally advantageous to find a so-called ``localized Wannier basis'' for the subspace spanned by these wavefunctions, i.e., a basis consisting of vectors whose support is concentrated on a small region of space. Such a basis exists under mild physical assumptions on the electronic system \cite{benzi_projector_decay}. Damle, Lin, and Ying \cite{damle_scdm} show how to find a localized Wannier basis by means of a CPQR factorization
\begin{equation}
\mPsi^*\mPi = \mQ\mR, \label{eqn:dft_cpqr}
\end{equation}
where $\mPsi \in \C^{n \times m}$ is the matrix of single-particle wavefunctions discretized onto a grid of size $n$.
\begin{figure}
    \centering
    \includegraphics[scale=.75]{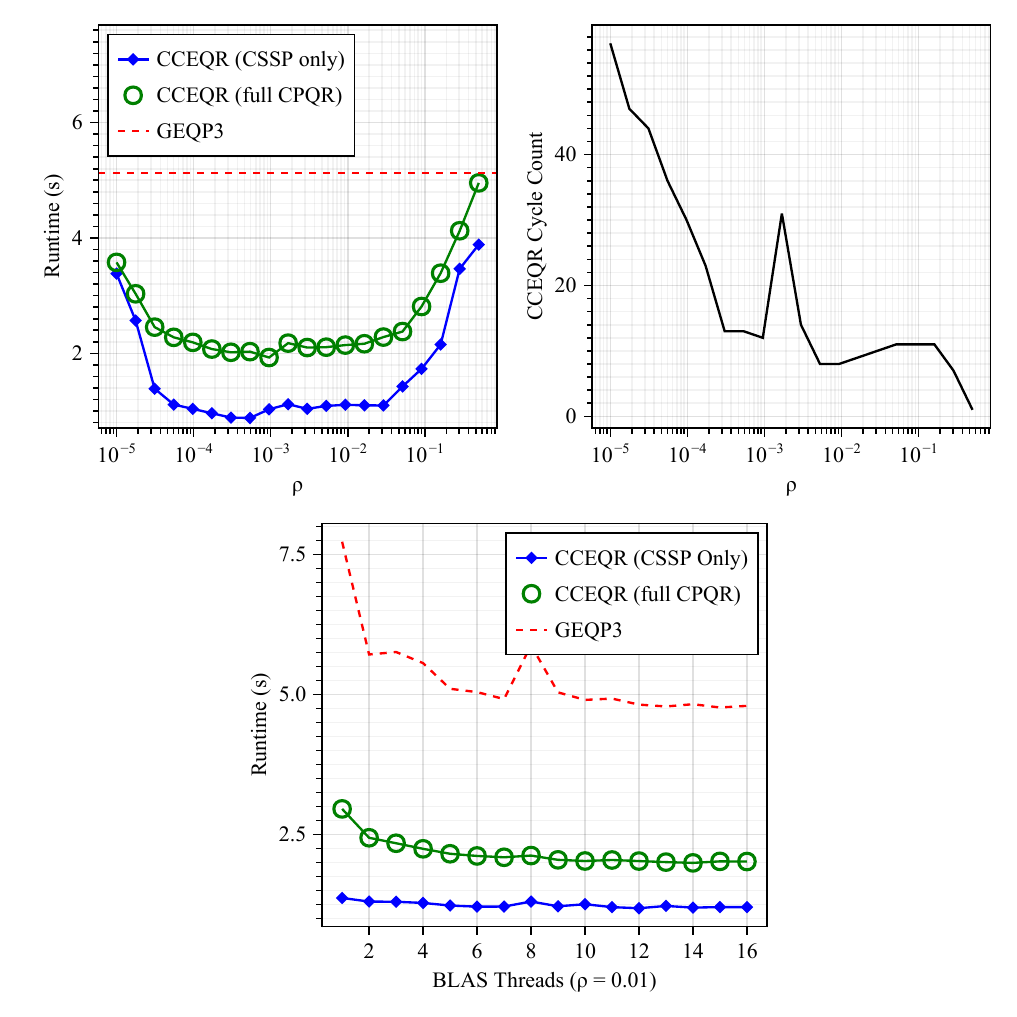}
    \caption{\underline{Top left}: median runtimes for CCEQR and GEQP3 over 30 trials on electronic wavefunctions for alkane with $m = 110$ and $n = 820{,}125$, selecting $k = m$ columns, with and without a full Householder reflection at the end of CCEQR. Runtime ratios for GEQP3 over CCEQR range from $1{.}32$ to $5{.}82$ (CSSP only) and from $1{.}03$ to $2{.}66$ (full CPQR); note that the full CPQR is not needed in this application. \underline{Top right}: cycle counts for CCEQR. \underline{Bottom}: median runtimes over 30 trials using different numbers of BLAS threads.} \label{fig:dft_alkane}
\end{figure}
\begin{figure}
    \centering
    \includegraphics[scale=.75]{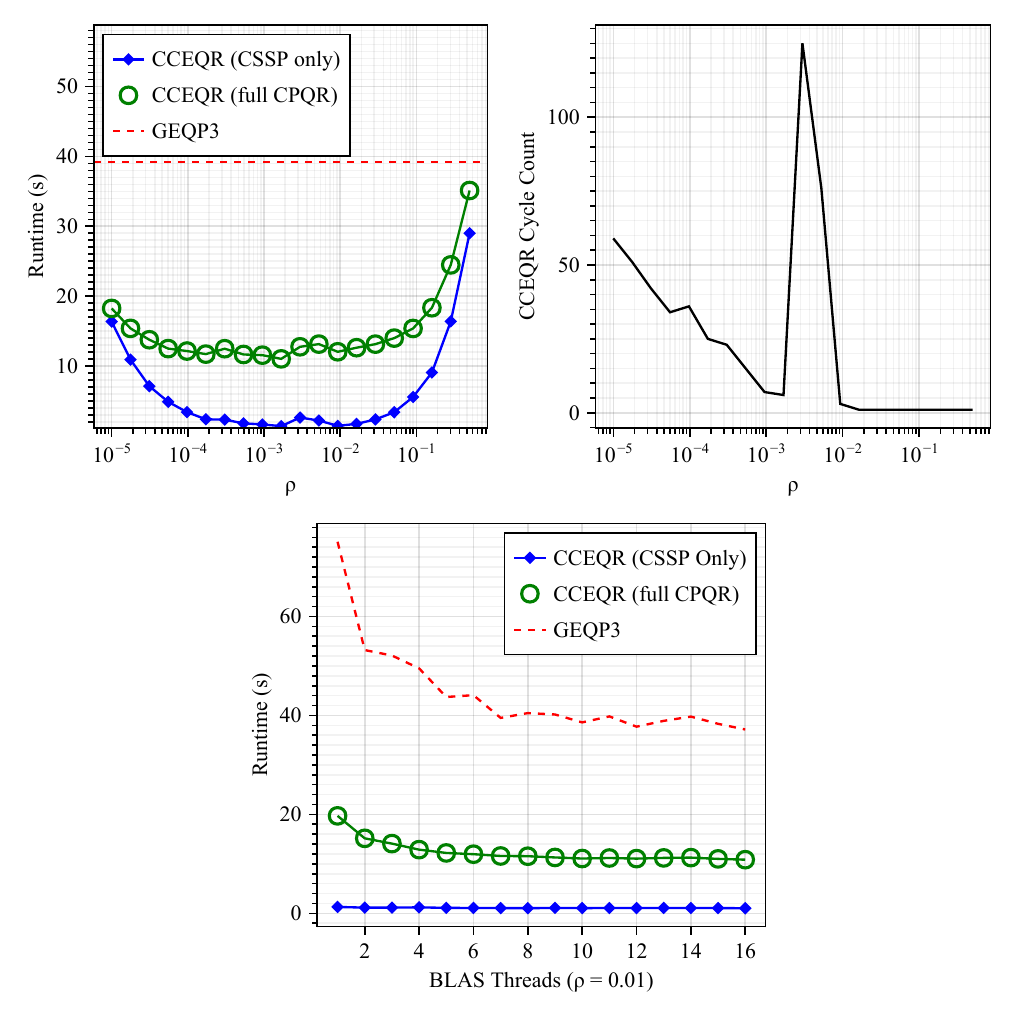}
    \caption{The same experiment as in \cref{fig:dft_alkane}, this time using 15 trials with wavefunctions from a water molecule with $m = 256$ and $n = 1{,}953{,}125$. See the main text for a discussion of the cycle count spike around $\rho \approx 3 \times 10^{-3}$. Runtime ratios for GEQP3 over CCEQR range from $1{.}35$ to $27{.}59$ (CSSP only) and from $1{.}12$ to $3{.}55$ (full CPQR). Note, as before, that the full CPQR is not necessary in this application.} \label{fig:dft_water}
\end{figure}

\Cref{fig:dft_alkane} compares the runtimes of GEQP3 and CCEQR for computing the column permutation in \cref{eqn:dft_cpqr}, where $\mPsi$ was generated from an alkane molecule with $m = 110$ and $n = 820{,}125$. For a broad range of values for $\rho$, both versions of CCEQR (with and without the final Householder reflection to produce all columns of $\mR$) perform significantly faster than GEQP3. Notably, this difference persists even when CCEQR takes many cycles to complete the permutation, indicating that the overhead of the ``expand'' step is not significant. Refer to \cref{subsec:clustering} for a more expansive discussion of how $\rho$ affects runtime. \Cref{fig:dft_alkane} also displays runtimes for varying numbers of BLAS threads. GEQP3 and the ``full CPQR'' variant of CCEQR see substantial benefits from increased multithreading, mainly because of accelerations in BLAS-3 Householder reflections applied to the full column set. For the ``CSSP only'' variant of CCEQR, which avoids applying reflections to all but a few columns, the benefits of multithreading are less significant.

\Cref{fig:dft_water} repeats this experiment for $\mPsi$ generated from a water molecule with $m = 256$ and $n = 1{,}953{,}125$. Here the performance difference between GEQP3 and CCEQR (without full Householder reflections) is even more pronounced, which should not be surprising given that $m$ and $n$ have both more-than-doubled relative to the alkane molecule. A striking feature in this experiment is the large spike in cycle count around $\rho \approx 3 \times 10^{-3}$. This disrupts our expectation that, as $\rho$ increases and the candidate blocks become larger, CCEQR should require \emph{fewer} cycles to find a complete column basis.

To understand this behavior, consider the state of CCEQR immediately after line \eqref{line:cceqr_collect_finalperm} in the ``collect'' step (\cref{alg:cceqr_collect}), where $b$ candidates have been selected out of $t$ tracked columns. The cycle-count spike in \cref{fig:dft_water} corresponds to a situation where there is a large cluster of non-candidate $j$'s for which $\vgamma(j)$ is very close to $\vgamma(s+1)$ (where $s+1$ is the index of the largest-residual candidate column). This can occur at the first cycle if the column norm distribution of the input matrix is near-uniform, and it can also occur more unpredictably at intermediate cycles if the residual matrix happens to acquire a near-uniform norm distribution.

In these cases, CCEQR is unlikely to commit any candidates into the skeleton \emph{except} for the one with largest residual. Indeed, the criterion for accepting a new skeleton column is that its squared-norm orthogonal to the skeleton exceeds the maximum of $\vgamma(j)$ over non-candidates; see equation \cref{eqn:cceqr_acceptance_number} and \cref{subsec:collect,subsec:commit}. In the situation described, given that the largest-residual candidate barely exceeds this maximum, orthogonalizing the remaining candidates against it is likely to decreased their residual norms below the required threshold. If a very large number of non-candidates are close to this maximum, then unless they are mostly colinear, CCEQR must work through a large number of cycles which each commit only a single column into the skeleton. Once enough of these cycles have taken place to sufficiently reduce residual norms in the cluster of large non-candidates, the algorithm can resume committing multiple columns at a time.

Notably, even when CCEQR encounters an ``obstacle'' of this sort, \cref{fig:dft_water} indicates that the overhead of excessive cycling does not significantly impact its runtime. This is not surprising, for if only a single candidate is committed then finishing the cycle requires only applying a \emph{single} Householder reflector to the tracked set.

\subsection{Gaussian Random Matrices} \label{subsec:random_gaussians}

\begin{figure}
    \centering
    \includegraphics[scale=.75]{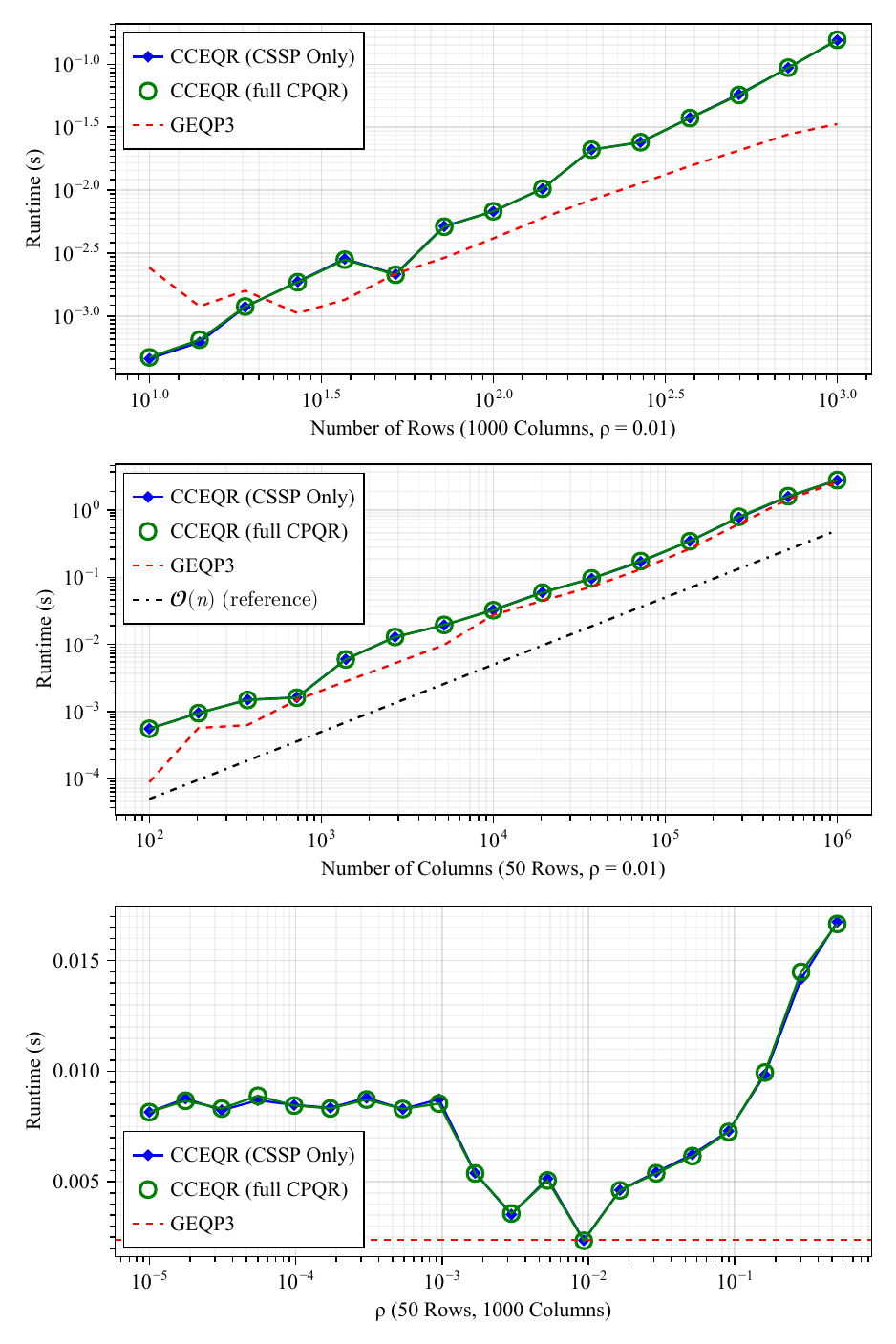}
    \caption{Median runtimes for CCEQR and GEQP3 over 30 trials on test matrices with i.i.d.\ standard Gaussian entries.} \label{fig:gaussians}
\end{figure}

Column subset selection problems do not always involve matrices whose column norm distribution is favorable to CCEQR. In light of this, \cref{fig:gaussians} compares CCEQR and GEQP3 on completely unstructured matrices whose entries are i.i.d.\ standard Gaussians. The results show, primarily, that CCEQR performs nearly the same (albeit slower) than GEQP3 on large unstructured problems. Importantly, CCEQR exhibits $\cO(n)$ runtime scaling for a fixed number of rows, just as GEQP3 does. Note that the runtime differences between CCEQR with and without full computation of $\mR$ are not significant. This is because for unstructured problems of this sort, CCEQR must bring nearly every column of the input matrix into the tracked set. In this case, the extra reflections needed to produce $\mR$ in full cost almost no extra work, since they are only applied to untracked columns at the end of the algorithm.
\begin{figure}
    \centering
    \includegraphics[scale=.75]{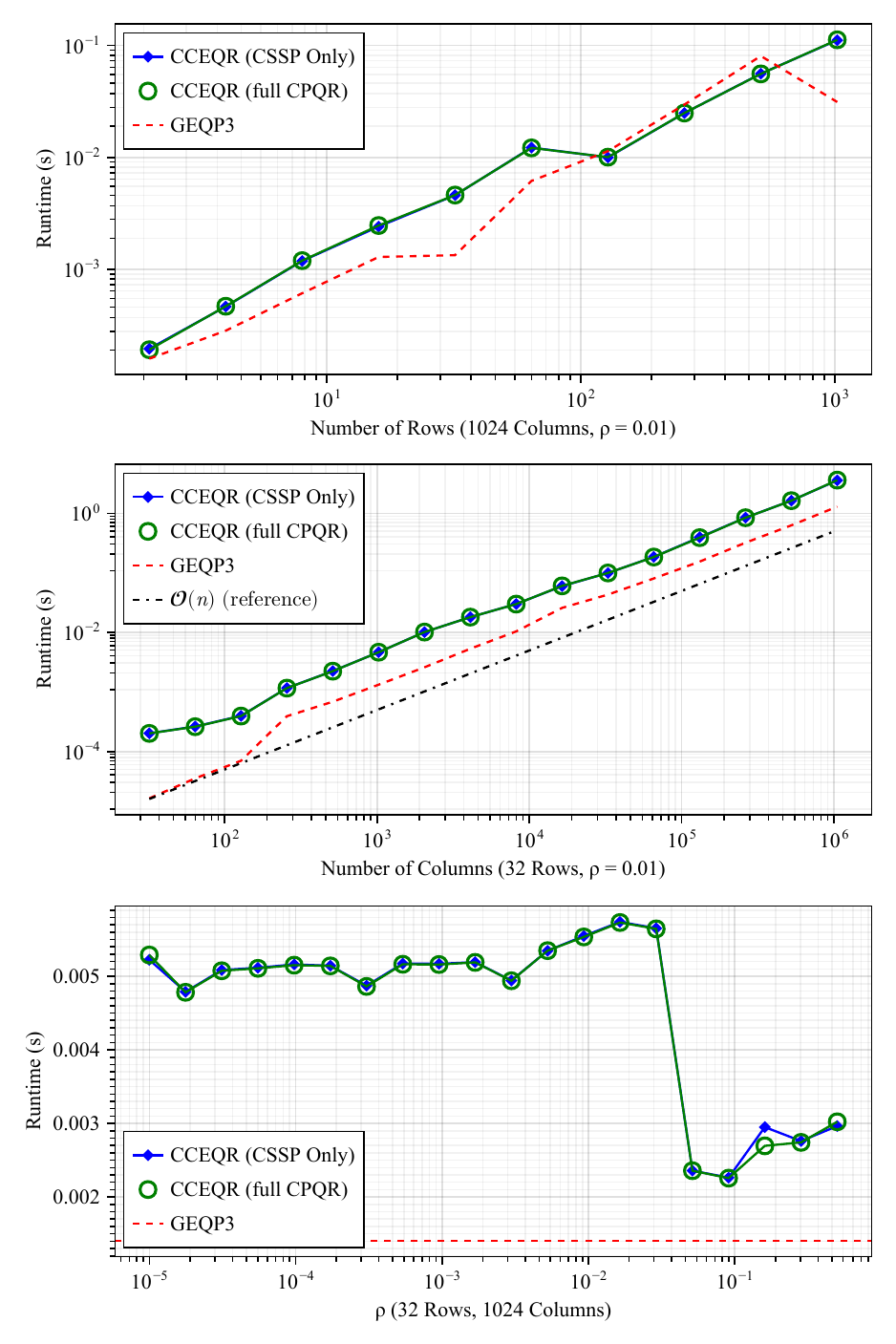}
    \caption{Median runtimes of CCEQR and GEQP3 over 30 trials on an ``adversarial'' $2^k \times 2^r$ Hadamard matrix, where $k = 5$ is fixed and $r$ is increased from 6 to 20.} \label{fig:hadamard}
\end{figure}

\subsection{Adversarial Hadamard Matrix} \label{subsec:hadamard_adversary}

Our final example is constructed ``adversarially'' to make CCEQR slower than GEQP3. This involves a matrix $\mH$ consisting of the first $2^k$ rows from a $2^r \times 2^r$ diadic Hadamard matrix. Such a matrix has two important properties:
\begin{enumerate}
\item all of its columns have equal norm, and
\item its columns can be partitioned into $2^k$ mutually orthogonal subsets of size $2^{r - k}$. Columns within a subset are all colinear.
\end{enumerate}
We have ordered the columns of $\mH$ so that every such subset appears in a contiguous block. To avoid issues related to column norm ties and floating point errors, we have also scaled the $j\nth$ column of $\mH$ by $1 + 1000(n - j + 1)\eta$ for $j = 1,\, \ldots,\, 2^r$, where $\eta = 2^{-52}$ is the machine epsilon.

This matrix is ``adversarial'' for CCEQR in the sense that, at any given cycle, there is an exceedingly large cluster of non-candidate columns whose residual squared norm almost matches that of the largest candidate column. Each cycle is therefore only able to commit a single column into the skeleton, for reasons discussed extensively in \cref{subsec:dft}. Our reordering of the columns into colinear blocks is meant to accentuate this behavior, since at every cycle, the orthogonalization step will disqualify remaining candidates in this block from being committed. Furthermore, because all column norms are essentially equal, CCEQR will need to bring every column into the tracked set, meaning significant work will be devoted to Householder reflections.

\Cref{fig:hadamard} shows runtimes of CCEQR and GEQP3 on this class of matrices, where $k$ varies from $1$ to $10$ and $r$ varies from 6 to 20. As in the random Gaussian test cases (cf.\ \cref{subsec:random_gaussians}), both CCEQR and GEQP3 exhibit $\cO(n)$ scaling for a fixed number of rows. Because the entire column set becomes tracked, there is also no significant runtime difference between CCEQR with and without full computation of $\mR$. It is encouraging that even for this ``adversarial'' example, CCEQR differs from GEQP3 in runtime by less than an order of magnitude. However, in contrast to the random Gaussian test cases, there is a persistent $\approx$ 2x runtime difference which is not eliminated by increasing the problem size. Is is possible that this gap could be made smaller by developing a more optimized implementation of CCEQR. But because this example is specifically designed to minimize the number of columns that can be committed at each cycle and maximize the tracked set size, thus creating the most inefficiency possible in CCEQR's control flow, we do not believe the runtime gap can be eliminated.
\section{Conclusions}
We have demonstrated an efficient CPQR-based column subset selection algorithm called CCEQR. This algorithm differs from existing CPQR-based rapid column selection algorithms in that (1) it is targeted toward matrices with far more columns than rows, (2) it is deterministic, and (3) it provably recovers the same column choice as the Golub Businger algorithm. Our algorithm is specifically designed for matrices whose column norm distribution has rapid decay. Using test matrices coming from applications in spectral clustering and density functional theory, which naturally have rapidly decaying column norms, we have demonstrated that CCEQR can run significantly faster than the LAPACK implementation of the Golub Businger algorithm (GEQP3). Although CCEQR can often outperform GEQP3 on the computation of a full column-pivoted QR factorization, the performance difference is most apparent for problems that only require computing the column permutation. We have also found that for problems whose column norm distribution is uniform, the performance difference between our algorithm and GEQP3 is small.

Future work could investigate whether or not a similar ``collect, commit, expand'' technique can be used to accelerate pivoted QR factorizations that are provably rank revealing, particularly the strong rank revealing QR factorization of Gu and Eisenstat \cite{gu_srrqr}. It may also be possible to further improve CCEQR's performance by implementing it in Fortran with direct BLAS and LAPACK calls, rather than the somewhat higher-level Julia code we have used in this paper.

\section{Acknowledgements}
RA and AD were partially supported by the National Science Foundation award DMS-2146079 and the Department of Energy Office of Science award DE-SC0025453. AD was also partially supported by the SciAI Center, funded by the Office of Naval Research under Grant Number N00014-23-1-2729.

\bibliographystyle{siamplain}
\bibliography{text/sources}

\appendix

\section{Proof of \cref{lemma:commit_rule}} \label{section:commit_rule_proof}
Let $\mQ,\, \mPi$ be the CPQR factors at the beginning of a given cycle of CCEQR, and define $\mR^{(0)} = \mQ^*\mA\mPi$, a matrix in $\mathrm{GB}(s)$ form. Let $\widehat{\mQ},\, \widehat{\mPi}$, and $c$ be defined as in \cref{lemma:commit_rule}. For the purposes of this proof, the important properties of $\widehat{\mQ}$ and $\widehat{\mPi}$ are as follows.
\begin{itemize}
\item Columns $(s+1) \mcol (s+t)$ of $\mR^{(0)}$ are arranged by $\widehat{\mPi}$ into an order determined using GEQP3 at line \eqref{line:collect_geqp3} of the ``collect'' step (\cref{alg:cceqr_collect}), and all other columns are left in place.
\item Rows $(s+1) \mcol m$ of $\mR^{(0)}$ are transformed by $\widehat{\mQ}$ using the first $c$ Householder reflectors from GEQP3, and all other rows are left unchanged.
\end{itemize}
Our goal in this section is to show that $\mR^{(1)} \defeq (\mQ\widehat{\mQ})^*\mA(\mPi\widehat{\mPi}) = \widehat{\mQ}^*\mR^{(0)}\widehat{\mPi}$ has $\mathrm{GB}(s+c)$ form. Before proceeding, let us check that the maximum defining $c$ in \cref{eqn:cceqr_acceptance_number} is over a nonempty set, so that $c \geq 1$ is well-defined. In the notation of \cref{subsec:collect}, let $\mB \defeq \mR^{(0)}((s+1) \mcol m,\, s + \vsigma(1 \mcol b))$ be the block of residual candidate columns, such that
\begin{equation*}
\widehat{\vp},\, \widehat{\vtau},\, \widehat{\mV},\, \widehat{\mR} \leftarrow \textsf{GEQP3}(\mB)
\end{equation*}
is the factorization performed at line \eqref{line:collect_geqp3} of the ``collect'' step (\cref{alg:cceqr_collect}). Recall that $c = \max\, \cI$, where
\begin{align*}
\cI &= \{ i \geq 1 \suchthat |\widehat{\mR}(i,i)|^2 \geq \delta \lor \mu \}, \\
\delta &= \max_{s+b < j \leq s+t} \| \mR^{(0)}\widehat{\mPi}((s+1) \mcol m,j) \|_2^2,\quad \mu = \max_{s+t < j \leq n} \| \mA\mPi(:,j) \|_2^2.
\end{align*}
Because the candidate block of the $b$ tracked columns with greatest residual norm (cf.\ \cref{subsec:collect}), we know that
\begin{align*}
|\widehat{\mR}(1,1)|^2 &= \max_{1 \leq j \leq b} \| \mB(:,j) \|_2^2 = \max_{s < j \leq s+t} \| \mR^{(0)}\widehat{\mPi}((s+1) \mcol m,j) \|_2^2 \geq \delta.
\end{align*}
Also, letting $\vgamma$ be the vector of residual column norms at the beginning of the given cycle, \cref{eqn:cceqr_norm_bounds} shows that
\begin{equation*}
|\widehat{\mR}(1,1)|^2 = \max_{s < j \leq s+t} \vgamma(j) \geq \max_{s+t < j \leq n} \| \mA\mPi(:,j) \|_2^2 = \mu.
\end{equation*}
We now see that $1 \in \cI$, so the maximum defining $c$ is well-defined.

Turning now to the main claim of \cref{lemma:commit_rule}, the first requirement of $\mathrm{GB}(s+c)$ form is that $\mR^{(1)}(:,1 \mcol (s+c))$ is upper-triangular. This follows from the construction of $\widehat{\mQ}$ and $\widehat{\mPi}$ with GEQP3. The second requirement is that for all indices $1 \leq i \leq s+c$ and $i \leq j \leq n$,
\begin{equation}
|\mR^{(1)}(i,i)| \geq \| \mR^{(1)}(i \mcol m,j) \|_2. \label{eqn:commit_rule_proof_target}
\end{equation}
We prove this by cases. 
\begin{enumerate}
    \item For $i \leq s$ we rely on the fact that $\mR^{(1)}(i,i) = \mR^{(0)}(i,i)$, and we consider two sub-cases.
    \begin{itemize}
        \item[1a.] If $j \leq s$ or $j > s+t$ then $\| \mR^{(1)}(i \mcol m,j) \|_2 = \| \mR^{(0)}(i \mcol m,j) \|_2$, as $\mR^{(1)}(i \mcol m,j)$ differs from $\mR^{(0)}(i \mcol m,j)$ by a rotation in the lower $m - s$ indices. This means that \cref{eqn:commit_rule_proof_target} follows from the $\mathrm{GB}(s)$ form of $\mR^{(0)}$.
        \item[1b.] If $s < j \leq s+t$, then there is a $p \in [s+1,\, s+t]$ determined by $\widehat{\mPi}$ such that $\mR^{(1)}(i \mcol m,j)$ and $\mR^{(0)}(i \mcol m,p)$ differ by only a rotation in the lower $m - s$ indices. Thus $\| \mR^{(1)}(i \mcol m,j) \|_2 = \| \mR^{(0)}(i \mcol m,p) \|_2$, so once again, \cref{eqn:commit_rule_proof_target} follows from the $\mathrm{GB}(s)$ form of $\mR^{(0)}$.
    \end{itemize}

    \item For $s < i \leq s+c$, we again consider sub-cases based on $j$.
    \begin{itemize}
        \item[2a.] If $j \leq s + b$ then
        \begin{align*}
        \mR^{(1)}(i,i) &= \widehat{\mR}(i-s,i-s) \\
        \mR^{(1)}(i \mcol m,j) &= \widehat{\mR}((i-s) \mcol (m-s), j-s)).
        \end{align*}
        In this case, \cref{eqn:commit_rule_proof_target} holds because GEQP3 outputs $\widehat{\mR}$ in $\mathrm{GB}(d)$ form (with $d = \min\{ m - s,\, b \}$).

        \item[2b.] If $s + b < j \leq s + t$ then, by \cref{eqn:cceqr_acceptance_number} and the $\mathrm{GB}(d)$ form of $\widehat{\mR}$,
        \begin{align*}
            |\mR^{(1)}(i,i)|^2 &= |\widehat{\mR}(i-s,i-s)|^2 \geq |\widehat{\mR}(c,c)|^2 \\
            &\geq \delta = \max_{s+b < j' \leq s+t} \| \mR^{(0)}\widehat{\mPi}((s+1) \mcol m,\, j') \|_2^2 \\
            &= \max_{s+b < j' \leq s+t} \| \mR^{(1)}((s+1) \mcol m,\, j') \|_2^2 \\
            &\geq \| \mR^{(1)}(i \mcol m,j) \|_2^2,
        \end{align*}
        as desired.

        \item[2c.] Finally, if $s+t < j$ then
        \begin{align*}
            |\mR^{(1)}(i,i)|^2 \geq |\widehat{\mR}(c,c)|^2 \geq \mu &= \max_{s+t < j' \leq n} \| \mA\mPi(:,j') \|_2^2 \\
            &= \max_{s+t < j' \leq n} \| \mR^{(1)}(:,j') \|_2^2 \\
            &\geq \| \mR^{(1)}(i \mcol m,j) \|_2^2,
        \end{align*}
        completing the proof.
    \end{itemize}
\end{enumerate}

\section{Efficient Updates of Householder Reflectors} \label{section:efficient_wy_details}
Given $m \times m$ unitary matrices $\mQ_1 = \mI - \mV_1\mT_1\mV_1^*$ and $\widehat{\mQ} = \mI - \widehat{\mQ}\widehat{\mT}\widehat{\mQ}^*$ in compact WY form, this section addresses the task of forming a compact WY representation of $\mQ_2 \defeq \mQ_1\widehat{\mQ}$ using BLAS-3 operations. In the context of CCEQR, $\mQ_1$ represents the unitary factor at the beginning of a given cycle, and $\widehat{\mQ}$ consists of the first $c$ Householder reflectors computed by GEQP3 in the ``collect'' stage (\cref{alg:cceqr_collect}). We write
\begin{align*}
\mV_1 = \begin{bmatrix} \vv_1 & \cdots & \vv_s \end{bmatrix},\quad \widehat{\mV} = \begin{bmatrix} \vv_{s+1} & \cdots & \vv_{s+c} \end{bmatrix},\quad \mV_2 = \begin{bmatrix} \vv_1 & \cdots & \vv_{s+c} \end{bmatrix}.
\end{align*}
By construction, $\mV_2$ is upper-triangular with unit diagonal. We have assumed here that $\widehat{\mV}$ is padded with zeros in the first few rows to have conformal dimensions with $\mV_1$, even though in \cref{subsec:commit}, $\widehat{\mV}$ denotes the raw output of GEQP3 without zero padding. \Cref{lemma:fast_wy_update} provides the needed update formulas.
\begin{lemma} \label{lemma:fast_wy_update}
If $\mI - \mV_2 \mT_2\mV_2^*$ is the compact WY form for $\mQ_2$, then
\begin{equation*}
\mT_2 = \begin{bmatrix}
\mT_1 & -\mT_1\mV_1^*\widehat{\mV}\widehat{\mT} \\
\mZero & \widehat{\mT}
\end{bmatrix}.
\end{equation*}
\end{lemma}

\begin{proof}
Inserting compact WY forms into $\mQ_2 = \mQ_1\widehat{\mQ}$, we have
\begin{equation*}
\mI - \mV_2\mT_2\mV_2^* = \mI - \mV_1\mT_1\mV_1^* - \widehat{\mV}\widehat{\mT}\widehat{\mV}^* + \mV_1\mT_1\mV_1^*\widehat{\mV}\widehat{\mT}\widehat{\mV}^*,
\end{equation*}
and because $\mV_2 = [\mV_1 \:\: \widehat{\mV}]$, this implies that
\begin{align*}
\mV_2\mT_2\mV_2^* &= \mV_1\mT_1\mV_1^* + \widehat{\mV}\widehat{\mT}\widehat{\mV}^* - \mV_1\mT_1\mV_1^*\widehat{\mV}\widehat{\mT}\widehat{\mV}^* \\
&= \mV_2 \begin{bmatrix}
\mT_1 & -\mT_1\mV_1^*\widehat{\mV}\widehat{\mT} \\
\mZero & \widehat{\mT}
\end{bmatrix} \mV_2^*.
\end{align*}
Because $\mV_2$ is a lower-triangular $m \times (s+c)$ matrix with unit diagonal, it has full column rank. We have used here the fact that, in the context of CCEQR, $s+c \leq k \leq m$. Therefore, multiplying by $\mV_2\pinv$ on the left and $(\mV_2^*)\pinv$ on the right proves the claim.
\end{proof}

\section{Numerical Results on Apple M2 Hardware} \label{section:apple_results} 
\begin{figure}
    \centering
    \includegraphics[scale=.8]{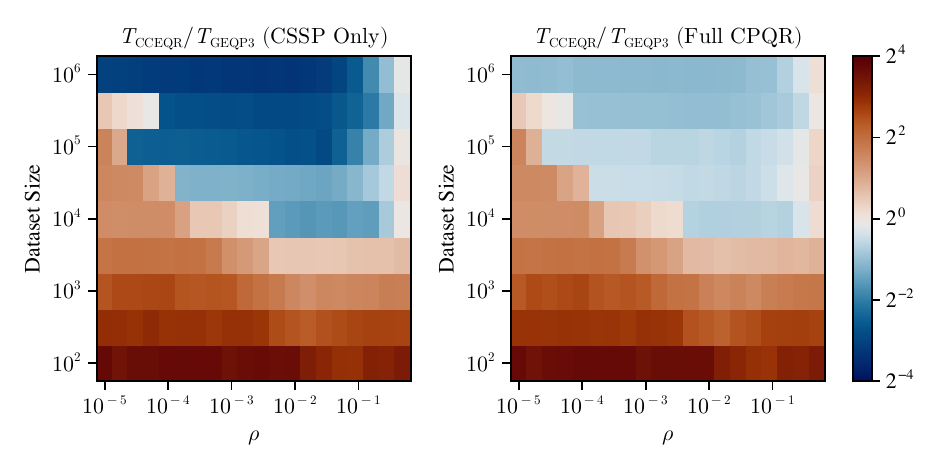}
    \caption{The same experiment as in \cref{fig:cluster_fixed_scale}, performed instead on an Apple M2 chip with the AppleAccelerate BLAS. In the left panel, plotted values range from $0{.}10$ to $13.69$. In the right panel, plottel values range from $0{.}47$ to $13{.}59$.} \label{fig:cluster_fixed_scale_apple}
\end{figure}

\begin{figure}
    \centering
    \includegraphics[scale=.8]{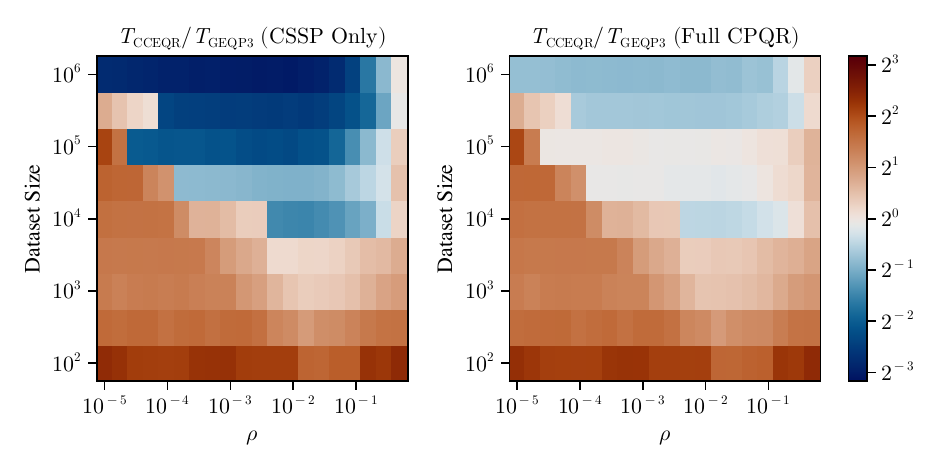}
    \caption{The same experiment as in \cref{fig:cluster_fixed_scale}, performed instead on an Apple M2 chip with OpenBLAS. In the left panel, plotted values range from $0{.}12$ to $5{.}43$. In the right panel, plottel values range from $0{.}56$ to $5{.}40$.} \label{fig:cluster_fixed_scale_openblas}
\end{figure}

For comparison with different hardware, we have repeated the experiment shown in \cref{fig:cluster_fixed_scale} on an Apple M2 chip using both the AppleAccelerate BLAS and OpenBLAS. The results are shown in \cref{fig:cluster_fixed_scale_apple,fig:cluster_fixed_scale_openblas}.

\end{document}